\documentclass[a4paper,11pt]{book}
\usepackage[english]{babel}
\usepackage{amsfonts}
\usepackage{amssymb}
\usepackage{amsmath}
\usepackage{latexsym}
\usepackage{mathrsfs}
\usepackage{mathtools}
\usepackage{rotating}
\usepackage{pdflscape}
\usepackage[utf8]{inputenc} 
\usepackage{textcomp} 
\usepackage{graphicx} 
\usepackage{epstopdf}
\usepackage{bm}  
\usepackage{pdfsync}  

\begin{document}
\centerline{\Large\bf Analytic theory of finite 
asymptotic expansions }
\par \centerline{\Large\bf in the real domain.} \vspace{5pt}
\par
\centerline{\Large\bf Part I:}\par
\centerline{\Large\bf  two-term expansions of differentiable functions.}

 \par\vspace{25pt}
\centerline{ANTONIO GRANATA}
 \par\vspace{20pt}

 \centerline{\footnotesize\sl Dipartimento di Matematica e Informatica, Universit\`a della
Calabria,} \par\centerline{\footnotesize\sl 87036 Rende (Cosenza), Italy, email:
antonio.granata@unical.it}

 \vspace{60pt}
 {\footnotesize {\bf Abstract. } It is our aim to establish a general analytic theory of asymptotic expansions of type 
$$
 f(x)=a_1\phi_1(x)+\dots+
a_n\phi_n(x)+o(\phi_n(x)),\ \ x\to x_0\,,\leqno(*)$$
where the given ordered $n$-tuple of real-valued
functions $(\phi_1\dots,\phi_n)$ forms an asymptotic scale at 
$x_o\in \overline {\mathbb R}$. 
By analytic theory, as opposed to the set of algebraic rules for manipulating finite 
asymptotic expansions, we mean sufficient and/or necessary conditions of general practical usefulness in order that $(*)$ hold true. Our theory is concerned with functions 
which are differentiable $(n-1)$ or $n$ times and the presented conditions involve
integro-differential operators acting on $f, \phi_1, \dots, \phi_n$. We essentially use
two approaches; one of them is based on canonical factorizations of $n$th-order disconjugate
 differential operators  and gives conditions expressed as convergence of certain improper integrals,
 very useful for applications. The other approach, valid for $(n-1)$-time differentiable functions starts
 from simple geometric considerations (as old as Newton's concept of limit tangent) and gives conditions expressed as the existence of finite limits, as $x\to x_{{}_0}$,of certain Wronskian determinants constructed with $f, \phi_1,\dots,\phi_n$. There is a link between the two approaches and it 
turns out that the integral conditions found via the factorizational approach have striking geometric meanings. Our theory extends to general expansions the theory of polynomial asymptotic expansions thoroughly investigated in a previous paper. In the first part of our work we study the case of two comparison functions $\phi_1,\phi_2$. The theoretical background for the two-term theory is much simpler than that for $n\ge3$ and, in addition, it is unavoidable to separate the treatments as the two-term formulas must be explicitly written lest they become unreadable.

The present e-paper coincides with the same-titled article published in 
``{\it Analysis Mathematica}, {\bf 37}(2011), 245-287,'' except for minor typographical alterations, for the addtion of a last section (\S9) discussing a nontrivial Tauberian problem, and for a list of corrections of misprints reported after the references, misprints that have been corrected in this e-paper.  }
  
\vspace{20pt}

 \vspace{10pt} {\bf AMS subject classifications.}
41A60, 34E05, 26C99.

\newpage
\centerline{\bf  Contents}
\vspace{10pt}
\begin{enumerate}
	\item Introduction
	\item Various approaches to the theory
	\item Basic assumptions and preliminary material
	\item The geometric approach
	\item The factorizational approach. Estimates of the remainder 
	\item The case of generalized convex functions 
	\item Proofs
	\item Example: the special case of powers
	\item Formal differentiation of a two-term asymptotic expansion:  a Tauberian result of  interpolatory type
\end{enumerate}

\vspace{15pt}
\centerline{\bf 1.\  Introduction}
\par\vspace{10pt}
Aim of our work is the establishing of a fairly complete theory of asymptotic expansions of type
$$
 f(x)=a_1\phi_1(x)+\dots+
a_n\phi_n(x)+o(\phi_n(x))\,, \quad x\to x_0\,,\quad n\geq 2\,,\leqno(1.1)
$$
where the comparison functions $\phi_i$ are supposed known in a neighborhood of $x_{{}_0}$ and forming an asymptotic scale at $x_{{}_0}$, i.e.
$$
\phi_1(x)>>\phi_2(x)>>\dots>>\phi_n(x)\,, \quad x\to x_0\,.\leqno(1.2)
$$
We deal with real-valued functions of one real variable. The simplest and first historical example of a relation (1.1)-(1.2), as all of us know, is Taylor's formula for which we have at disposal not only the elementary rules for manipulating the corresponding asymptotic relations but also the various forms of Taylor's theorem which give simple analytic conditions on $f$ sufficient for the validity of Taylor's formula of a certain order either with a simple asymptotic estimate or with some representation of the remainder. As far as general asymptotic expansions are concerned the current state of affairs is that we are able to perform practical manipulations, see e.g. Dieudonn\'e [2; ch. III], but no analogue of Taylor's formula is explicitly stated in the literature except in the case of expansions in real powers, i.e. $\phi_i(x)=x^{\alpha_i}$. In [5] the author collected and systematized various scattered results concerning polynomial expansions
$$
 f(x)=a_n x^n+\dots+
a_1 x+a_0+o(1)\,, \quad x\to +\infty\,,\leqno(1.3)
$$
with an eye to highlight the geometric approach and to link different approaches by a unique thread. In [6; 7] the author developed a theory for expansions in real powers
$$
f(x)=a_1 x^{\alpha_1}+\dots+a_n x^{\alpha_n}+o(x^{\alpha_n})\,, \quad x\to +\infty\,;\quad\alpha_1>\dots>\alpha_n\,,\leqno(1.4)
$$
with the aim of obtaining complete and applicable results about the formal dif\-fe\-rentiation of (1.4), results not obtainable by any of the classical approaches used for formal dif\-fe\-ren\-tia\-tion of the asymptotic relations $f(x)=O(x^{\gamma})$ or $f(x)=o(x^{\gamma})$.
\par
Now we intend to develop a complete theory of expansions (1.1)-(1.2) for differentiable functions.Our title 
``{\sl Analitic theory of ... }" is meant as opposed to the ``{\sl algebraic theory }" i.e. the set of rules for manipulating finite asymptotic expansions. Our work deals with functions which are differentiable a sufficient number of times and the exhibited sufficient and/or necessary conditions involve certain dif\-fe\-ren\-tial operators in the same way as Taylor's formula involves $n$th-order derivatives. Different approaches are used and complementary results are obtained but there is one guiding thread: the theory of P\'olya-Mammana fa\-cto\-ri\-za\-tions of linear ordinary differential operators in its latest developments, due to Trench [13] and the author [3; 4], concerning canonical factorizations. The $n$-tuple $(\phi_1,\dots,\phi_n)$ is subject to the practically mild restriction of forming a Chebyshev system on a one-sided neigh\-borhood of $x_{{}_0}$ and this yields ``natural" generalizations of Taylor's formula. Adapting on old method by Haupt [10] we find a geometric characterization of a certain asymptotic situation involving both (1.1) and suitable relations for the derivatives of $f$.
\par
Our exposition is split in two parts, according as $n=2$ or $n\geq 3$, for technical and practical reasons. In fact the theory for $n\geq 3$ requires the results for $n=2$ as some proofs are by induction on $n$; moreover statements and formulas for $n=2$ must be written out explicitly to avoid misinterpretations of the more complicated higher-order formulas.
\par
Propositions are numbered consecutively irrespective of their labelling as theorem, lemma etc..
\par\vspace{20pt}
\centerline{\bf 2.\ Various approaches to the theory}
\par\vspace{10pt}
For a general two-term expansion
$$
 f(x)=a_1\phi_1(x)+a_2\phi_2(x)+o(\phi_2(x))\,, \quad x\to x_0\,,\leqno(2.1)
$$
where $\phi_1, \phi_2$ do not vanish on a suitable deleted neighborhood of $x_{{}_0}$ and $\phi_1(x)>>\phi_2(x)$, $x\to x_{{}_0}$, we have the characterizing relations
$$
a_1=\lim_{x\to x_0} f(x)/\phi_1(x)\,;\quad a_2=\lim_{x\to x_0}\left.[f(x)-a_1/\phi_1(x)]\right/\phi_2(x)\,,\leqno(2.2)
$$
and we know that they can seldom be used in meaningful applications other than elementary cases. For this reason we look for sufficient, and possibly necessary, conditions of a quite different nature for (2.1) to hold. We point out four different approaches which are fit for the $n$-term theory as well.

\textbf{1.\ The naive approach.}\ If we try to apply L'Hospital's rule to evaluate the second limit (2.2), writing
$$
(f-a_1\phi_1)/\phi_2\equiv {f/\phi_1-a_1\over (\phi_2/\phi_1)}\,,
$$
we see that a sufficient condition for (2.2) to hold is the existence as finite numbers of the first limit (2.2) and of 
$$
a_2:=\lim_{x\to x_0} {(f/\phi_1)'\over (\phi_2/\phi_1)'}\,.\leqno(2.3)
$$\par
We label this approach as ``naive" because  its elementary idea leads us to replace the second limit (2.2) by a seemingly more complicated limit: maybe this approach is a blind alley and we should better try other paths. However it is elemenary to prove that the existence of the limit (2.3) implies the existence of the first limit (2.2) due to condition $\phi_2(x)/\phi_1(x)=o(1)$, $x\to x_{{}_0}$, and to condition ``$(\phi_2(x)/\phi_1(x))'$ strictly one-signed  on a neighborhood of $x_{{}_0}$'', which is necessary for the application of L'Hospital's rule. Hence the pair of conditions (2.2) is implied by the single condition (2.3) defining $a_2$ independently from $a_1$, at least under the additional restriction on the sign of $(\phi_2/\phi_1)'$. So far our investigation is nothing but an elementary exercise but we shall point out the asymptotic, the geometric and the analytic meaning of (2.3).

\textbf{2.\ The algebraic approach: formal differentiation of asymptotic expansions.}\quad Let us examine the case wherein the remainder in (2.1) is identically zero i.e. the given function $f$ coincides, at least in a neighborhood of $x_{{}_0}$, with a ``generalized polynomial"
$$
f(x)\equiv a_1\phi_1(x)+a_2\phi_2(x)\,.\leqno(2.4)
$$
\par Besides $a_1$, defined by the first relation (2.2), we may express $a_2$ independently from $a_1$ through a suitable differential operator. In fact we get in sequence from (2.4):
$$
f(x)/\phi_1(x)=a_1+a_2(\phi_2(x)/\phi_1(x))\,,\leqno(2.5)
$$
$$
(f(x)/\phi_1(x))'=a_2(\phi_2(x)/\phi_1(x))'\,,\leqno(2.6)
$$
$$
a_2\equiv {(f(x)/\phi_1(x))'\over (\phi_2(x)/\phi_1(x))'}\,,\leqno(2.7)
$$
provided the last expression takes a meaning on some interval. We now apply the same procedure to (2.1) first obtaining 
$$
f(x)/\phi_1(x)=a_1+a_2(\phi_2(x)/\phi_1(x))+o(\phi_2(x)/\phi_1(x))\,,\quad x\to x_0\,,\leqno(2.8)
$$
and then we conjecture that formal differentiation both sides of (2.8) may yield
$$
(f(x)/\phi_1(x))'=a_2(\phi_2(x)/\phi_1(x))'+o\left((\phi_2(x)/\phi_1(x))'\right)\,,\quad x\to x_0\,,\leqno(2.9)
$$
at least under ``reasonable conditions". Relation (2.9) is an equivalent formulation of (2.3) but our procedure leads us to interpret (2.9) as a relation obtained from (2.1) by formal application of a certain first-order differential operator. We have already mentioned that (2.3) implies (2.2); our present approach suggests other natural questions: 
\par (i) Does the existence of the limit (2.3) characterize the pair (2.8)-(2-9), i.e. the pair (2.1)-(2.9)?
\par (ii) Does this contingency occur in meaningful situations or does it occur in exceptional situations only?
\par (iii) What about the limit
$$
\lim_{x\to x_0} {(f(x)/\phi_2(x))'\over (\phi_1(x)/\phi_2(x))'}\leqno(2.10)
$$
whimsically obtained from (2.3) by interchanging the roles of $\phi_1$, $\phi_2$? Is it an unnatural quantity to be taken into consideration or has it a meaning in our context?
\par
All these questions will receive answers in this paper. Once again notice that the present approach considers the expansion (2.1) not by itself but matched to (2.9), which is obtained by a certain process of formal differentiation.

\textbf{3.\ The analytic or factorizational approach.}\ The idea is to use integro-\--differential representations of $f$ from which one can easily infer sufficient and/or necessary conditions for (2.1) to hold in the very same way as an expansion of any of the simple types
$$
f(x)=a_1+a_2x+o(x)\,,\; x\to 0\,;\quad f(x)=a_1x+a_2+o(1)\,,\; x\to +\infty\,,\leqno(2.11)
$$
can be studied starting either from the familiar representation
$$
f(x)=c_1x+c_2+\int\limits_T^x dt \int\limits_T^t f'' (\tau) d\tau\leqno(2.12)
$$
or from the less usual representation
$$
f(x)=\overline c_1x+\overline c_2+x \int\limits_T^x t^{-2} dt \int\limits_T^t \tau f'' (\tau) d\tau\,.\leqno(2.13)
$$
\par
The procedure goes as follows. Consider a second-order linear ordinary differential operator $L$ whose null-set coincides with \textsl{span} $(\phi_1, \phi_2)$ and which can be explicitly inverted. The most favourable circumstance is whenever $L$ admits of a factorization 
$$
L[u]\equiv p_2[p_1(p_0 u)']'\leqno(2.14)
$$
with suitable nowhere-vanishing functions $p_i$; in this case one can write down a corresponding integral representation of $f$ in terms of $L[f]$:
$$
f(x)=c_1\phi_1(x)+c_2\phi_2(x)+{1\over p_0(x)}\int\limits_T^x {dt_1\over p_1(t_1)}\int\limits_T^{t_1}{L[f(t_2)]\over p_2(t_2)}dt_2\,.\leqno(2.15)
$$
\par
Now one tries to obtain analytic characterizations of (2.1),or of the pair (2.1)-(2.9), or of other pairs of asymptotic relations via integrability conditions on $L[f]$.  This approach automatically gives integral representations of the remainders. An exhaustive investigation involves the two types of canonical factorizations available for $L$ and described in \S 3.

\textbf{4.\ The geometric approach.}\ The idea is nothing but Newton's concept of ``limit tangent" to the graph of a function as the point of contact goes to infinity. The straight line tangent to the graph of $f$ at a generic point $(t, f(t))$ is represented by equation
$$
y=f(t)+f'(t)(x-t)\equiv f'(t)x+[f(t)-tf'(t)]\,,\leqno(2.16)
$$
and its limit position as $t\to +\infty$ may be defined as the line $y=a_1x+a_2$ where
$$
a_1:=\lim_{t\to +\infty} f'(t)\;\,;\; a_2:=\lim_{t\to +\infty} [f(t)-tf'(t)]\,.\leqno(2.17)
$$
\par
A satisfying treatment goes back to the German geometer Haupt [8] almost one century ago. Applying the same idea to (2.1) one first chooses, among all linear combinations $c_1\phi_1(x)+c_2\phi_2(x)$, that special one which has a first-order contact with $f$ at a generic pont $t$ and which is characterized by certain coefficients $c_1=f^*_1(t)$, $c_2=f^*_2(t)$. Then one investigates relationships between the existence, as finite numbers, of any one or both of the limits $\lim_{t\to x_o} f^*_i(t)$, $i=1,2$, and the validity of the expansion (2.1). Such a procedure provides analytic characterizations of various asymptotic situations by means of geometrically-meaningful quantities; the $\lim_{t\to x_o} f^*_2(t)$ turns out to be more meaningful than the $\lim_{t\to x_o} f^*_1(t)$. The elementary case $(\phi_1(x),\phi_2(x))\equiv(1,x)$ as well as the
 general polynomial case $(\phi_1(x),\dots,\phi_n(x))\equiv(1,x,\dots, x^{n-1})$, studied in [5], suggest the introduction of another meaningful geometric quantity $F^*(t)$, definition 4.1 below, which represents the ordinate of the intersection point between the graph of the osculating curve $y=f^*_1(t)\phi_1(x)+f^*_2(t)\phi_2(x)$ at the generic point $(t, f(t))$ and a fixed vertical line $x=T$. The main result states the equivalence of the existence as finite quantities of any one of the following three limits:
$$
\lim_{t\to x_0} F^*(t)\;\,;\; \lim_{t\to x_0} f_2^*(t)\;\,;\;\lim_{t\to x_0} (f_1^*(t), f^*_2(t))\,.\leqno(2.18)
$$
\par
This contingency defines the ``limit" of the osculating curve and is also equivalent to the pair (2.1)-(2.9).
\par
\centerline{*\quad *\quad *}
\par
Following the third and fourth approaches we shall give substance to the first and second approaches. The quantity $(f(t)/\phi_1(t))'/(\phi_1(t)/\phi_2(t))'$, found in the naive approach, turns out to be the coefficient $f_2^*(t)$ whereas the ``whimsical'' quantity appearing in the limit (2.10)
is the coefficient $f^*_1(t)$: see \S 4. Moreover, if all the involved functions are supposed to be of class $AC^1(I)$, i.e. endowed with first-order derivatives absolutely continuous on a neighborhood $I$ of $x_{{}_0}$, then $(f^*_i)'(t)=q_i(t)\cdot L[f(t)]$, $i=1,2$, where $q_i$ are suitable nowhere-vanishing functions and $L$ is a second-order linear ordinary differential operator such that $ker\, L= span (\phi_1, \phi_2)$. Expressing $f^*_i$ as integral operators acting on $L[f]$ we transmute all the geometric conditions into simpler and practically useful analytic conditions. Last, but not least, a unique theoretical link is found for all the approaches, namely the theory of canonical factorizations of second-order differential operators.
\par
The theory we shall develop for two-term expansions contains all the essential ideas for the $n$-term expansions. In the second part of our work we shall point out the contributions of other authors to the $n$th-order theory, especially Kudryavtsev's Lagrangian approach (which is nothing but the geometric approach) to a larger class of expansions.

\vspace{20pt}\centerline{\textbf{3.\ Basic assumptions and preliminary material}}

\vspace{10pt}
In studying (2.1) the real-valued functions $\phi_1, \phi_2, f$ are supposed to be defined on a deleted one-sided neighborhood $I$ of $x_{{}_0}$ and, for definiteness, we suppose that $I$ is left-sided and $x_{{}_0}\leq +\infty$. Every limit process ``$\lim_{x\to x_{{}_0}}\cdots$" stands for ``$\lim_{x\to x_{{}_0}, x\in I}\cdots$". If $I$ is any interval in $\mathbb{R}$ the symbol $AC^k(I)$ denotes the space of all $f$'s such that $f\in C^k(I)$ and $f^{(k)}$ is absolutely continuous on every compact interval of $I$: i.e. $f\in AC^k(I)\Leftrightarrow f^{(k)}\in AC (I)$. $AC^\circ (I)\equiv AC(I)$. Whenever $f\in AC(I)$, writing ``$\lim_{x\to x_{{}_0}} f'(x)$" stands for ``$\lim_{x\to x_{{}_0}, x\in \widetilde I} f'(x)$" where $\widetilde I$ is the subset of $I$ where $f'$ exists as a finite number; when applying L'Hospital's rule in such a context we always use Ostrowski's version [11] valid for absolutely conti\-nuous functions. The symbols $f\in L^1 (I)$, $f\in L^1_{\hbox{\scriptsize{loc}}}(I)$ respectively denote that $f$ is Lebesgue-summable on $I$ or on every compact interval of $I$; $f$ \textsl{integrable} on $[T, x_{{}_0}[$ means that $f \in L^1_{\hbox{\scriptsize{loc}}}[T, x_{{}_0}[$ and the improper integral
$$
\int\limits_T^{x_0} f \equiv \int\limits_T^{\to x_0} f:=\lim_{x\to x_0} \int\limits_T^x f\quad\hbox{converges}\,.
$$
$\overline{\mathbb R}$ is the extended real line, $\overline{\mathbb R} :={\mathbb R} \cup\{\pm \infty\}$.
\par\textbf{Basic assumptions on } $(\phi_1, \phi_2)$:
$$
\phi_1, \phi_2\in C^1(I)\;\,;\; I:=[T, x_0[\,,\quad T\in \mathbb{R}\,;\leqno(3.1)_1
$$
$$
\phi_1(x) >> \phi_2(x)\,, x\to x_0\quad \hbox{i.e. }\quad \phi_2(x)=o(\phi_1(x))\,,\; x\to x_0\,;\leqno(3.1)_2
$$
$$
\phi_1(x)\,,\; \phi_2(x)\neq 0\quad \forall x\in I\,;\leqno(3.1)_3
$$
$$
W(x):=W(\phi_1(x), \phi_2(x))\equiv W(\phi_1, \phi_2; x)\neq 0\quad \forall x\in I\,,\leqno(3.1)_4
$$
\textsl{where} $W(\phi_1, \phi_2)$ \textsl{is the Wronskian determinant of} $\phi_1, \phi_2$.
\par Immediate consequences of assumptions (3.1) are that $\phi_1, \phi_2$ are linearly independent on $I$; that $\phi_1, \phi_2, W$ are strictly one-signed on $I$ and that
$$
\int\limits^{x_0} (\phi_2(t)/\phi_1(t))' dt \quad converges\,,\leqno(3.2)
$$
$$
\phi_2(x)=-\phi_1(x)\cdot \int\limits^{x_0}_x (\phi_2(t)/\phi_1(t))' dt\;\,,\;\; x\in I\,.\leqno(3.3)
$$
\par
It is also known that any linear combination of $\phi_1, \phi_2$ is either $\equiv 0$ on $I$ or has at most one zero on $I$: see, e.g., Coppel [1; prop. 5, p. 89].
\par
It will be specified in \S 6 that, by changing the signs of $\phi_1, \phi_2$ if necessary, the three conditions (3.1)$_1$, (3.1)$_3$, (3.1)$_4$, state that the ordered pair of functions $(\phi_1, \phi_2)$ is a Chebyshev system on $I$ and this is the theoretical framework of our theory making possible the geometric and the analytic approaches and not only the mechanical application of L'Hospital's rule to evaluate the second limit (2.2).
\par\vspace{5pt}\textbf{Strenghtened basic assumptions on $(\phi_1, \phi_2)$:}\par
\textsl{All assumptions } (3.1) \textsl{plus the stronger regularity condition}
$$
\phi_1, \phi_2 \in AC^1 (I)\leqno(3.4)
$$
\par 
In this case, besides the above-mentioned consequences, there exists a unique second-order linear  ordinary differential operator $L$
$$
L[u]:= u''+a_1(x)u'+a_2(x)u\;\,;\; a_i\in L^1_{\scriptstyle{loc}}(I)\leqno(3.5)
$$
such that
$$
ker\ L=span\ (\phi_1, \phi_2)\,.\leqno(3.6)
$$
\par
If $\phi_1,\phi_2\in C^2(I)$ then $a_1, a_2\in C^\circ (I)$. Condition (3.1)$_3$ now implies that the operator $L$ is disconjugate on $I$, Coppel [1; th. 1, p. 5], hence it admits of a P\'olya-Mammana factorization on $I$, i.e.
$$
L[u]\equiv p_2[p_1(p_0 u)']'\qquad \forall u\in AC^1(I)\,,\leqno(3.7)
$$
where $p_0, p_1, p_2$ are suitable functions strictly positive on $I$ and satisfying the regularity conditions:
$$
p_0\in AC^1(I)\;\,;\quad p_1, p_2\in AC^0(I)\,.\leqno(3.8)
$$
\par
For a given factorization (3.7) any function $f\in AC^1(I)$ admits of an integral representation of type
$$
f(x)=c_1\phi_1(x)+c_2\phi_2(x)+{1\over p_0(x)}\int\limits_{T_0}^x{dt_1\over p_1 (t_1)}\int\limits_{T_0}^{t_1}{L[f(t_2)]\over p_2 (t_2)} dt_2\,,\; x\in I\,,\leqno(3.9)
$$
where $T_0$ is arbitrarily chosen in $I$ and $c_1, c_2$ are suitable constants depending on $f$, $\phi_i, p_i, T_0$.
\par
\textbf{Factorizations of the operator $L$.}\ \ 
Following the terminology used in [3], factorization (3.7) is called a ``\textsl{canonical factorization of type (I) or of type (II), at} $x_{{}_0}$" according as the improper integral $\int^{x_{{}_0}} (1/p_1)$  respectively diverges or converges.
\par\vspace{5pt}\textbf{Lemma 3.1}\quad (see [4; th. 6.1]. \textsl{For each solution} $\phi$ \textsl{to} $L[u]=0$ \textsl{strictly positive on a left-sided neighborhood of} $x_{{}_0}$, $J\subset [T, x_{{}_0}[$, \textsl{there exists an} ``\textsl{essentially unique}" \textsl{factorization} (3.7)-(3.8) \textsl{on} $J$ \textsl{such that} $p_0=1/\phi$, \textsl{namely}
$$
L[u]\equiv {W(x)\over \phi(x)}\left[{(\phi(x))^2\over W(x)}\left({u\over \phi(x)}\right)'\right]'\quad \forall\ u\in AC^1(J)\,.\leqno(3.10)
$$
\par
(Here ``essentially unique"' means that the functions $p_i$ are determined up to constant factors).
\par\vspace{5pt}\textbf{Lemma 3.2}\quad \textsl{Factorization} (3.10) \textsl{is a canonical factorization of type (I) at} $x_{{}_0}$  \textsl{iff}
$$
\phi(x)\sim c\phi_2(x)\,,\quad x\to x_0\quad (c\neq 0)\,;\leqno(3.11)
$$
\textsl{and it is a canonical factorization of type (II) at} $x_{{}_0}$ \textsl{iff}
$$
\phi(x)\sim c\phi_1(x)\,,\quad x\to x_0\quad (c\neq 0)\,,\leqno(3.12)
$$
\textsl{where} $c$ \textsl{denotes a suitable constant. Contingency} (3.11) \textsl{occurs iff} $\phi(x)\equiv c \phi_2(x)$ \textsl{on} $I$; \textsl{hence on a fixed left-sided neighborhood of} $x_{{}_0}$ \textsl{there exists only} ``\textsl{one}" \textsl{canonical factorization of type (I) at} $x_{{}_0}$.
\par\vspace{5pt}
\textbf{Proof.}\quad By (3.1)$_2$ the two circumstances (3.11) and (3.12) are the only ones that can occur as far as the asymptotic behavior of $\phi(x)$ at $x_{{}_0}$ is concerned. If $\phi\sim c\phi_1$ then
$$
{W(x)\over (\phi(x))^2}\sim c^{-2}{W(\phi_1(x), \phi_2(x))\over (\phi_1(x))^2}\equiv c^{-2}(\phi_2(x)/\phi_1(x))'\,,\quad x\to x_0\,,\leqno(3.13)
$$
and this implies, by (3.2), the convergence of the improper integral
$$
\int\limits^{x_0} W(t) (\phi(t))^{-2}dt\,.\leqno(3.14)
$$
\par
If $\phi\sim c\phi_2$ then (3.1)$_2$ implies that $\phi(x)=c\phi_2(x)$ on $J$ and relation (3.13) is replaced by
$$
{W(x)\over (\phi(x))^2}=-c^{-2}(\phi_1(x)/\phi_2(x))'\,,\quad x\to x_0\,,\leqno(3.15)
$$
which implies the divergence of the integral (3.14) as $\lim_{x\to x_{{}_0}}\phi_1(x)/\phi_2(x)=\pm\infty$.
\par\hfill{$\Box$}
\par\vspace{5pt}
\par
We bring to the reader's attention the fact that for a generic operator $L$ disconjugate on an open interval $]T, x_{{}_0}[$ it may happen that there exist no solution $\phi$ to $L[u]=0$ strictly positive on the whole interval $]T, x_{{}_0}[$ and satisfying (3.12): just think of the operator $u''$ on $(-\infty, +\infty)$. However this is granted if we restrict $]T, x_{{}_0}[$ to $]T+\epsilon, x_{{}_0}[$ or $]T, x_{{}_0}-\epsilon[$, with an arbitrary $\epsilon >0$, as implied by [3; th.2.2, p. 162]. Our assumption (3.1)$_3$ only means that we are restricting the original interval if necessary. We shall develop our theory taking advantage of both types of canonical factorizations. It is immaterial whether in a factorization, either in this paper or in practical applications, some coefficient $p_i$ happens to be strictly negative.
\par
\centerline{*\quad *\quad *}
\par
The following trivial formulas are reported for the sole purpose of later references. Any of the following three notations will be used to denote the Wronskian of two functions
$$
W(f,g)\equiv W(f,g; x)\equiv W(f(x), g(x)):= f(x) g'(x)-f'(x) g(x)\,;\leqno(3.16)
$$
$$
(f/g)'=-(g/f)'\cdot (f/g)^2\quad\hbox{ if }\; f,g\neq 0\,;\leqno(3.17)
$$
$$
W(f,g)= f^2\cdot (g/f)'\quad\hbox{ if }\;f\neq 0\,;\quad W(f,g)=-g^2\cdot (f/g)'\quad\hbox{ if }\;g\neq 0\,; 
\leqno(3.18)
$$
$$
{d\over dx} W(f,g; x)=f(x) g''(x)-f''(x) g(x)\,;\leqno(3.19)
$$
$$
\begin{array}{ll}
\kern-1.3cm(3.20)\qquad &W(a_1f_1+a_2f_2, b_1g_1+b_2g_2)= a_1b_1 W(f_1,g_1)+a_1b_2 W(f_1,g_2)+\\ \\
&+a_2b_1 W(f_2,g_1)+a_2b_2 W(f_2,g_2)\,,\qquad (a_i\,, b_i=\hbox{constants})\,.\\
\end{array}$$

If $f_1, f_2, g_1, g_2$ are linked by the relations
$$\begin{cases}
f_1(x)=a_{11} g_1(x)+a_{12} g_2(x))\cr\cr
f_2(x)=a_{21} g_1(x)+a_{22} g_2(x))\cr \qquad (a_{ij}=\hbox{constants})
\end{cases}\leqno(3.21)$$
then
$$
W(f_1(x),f_2(x))=
\left\vert\begin{array}{cc}
a_{11} & a_{12}\\
a_{21} & a_{22}\\
\end{array}\right\vert
\cdot W(g_1(x), g_2(x))\,.\leqno(3.22)
$$
\par
According to (3.1)$_4$ the symbol $W(x)$ always refer to the Wronskian of the comparison functions $\phi_1, \phi_2$ fixed in any particular context.
\par\vspace{20pt}\centerline{\textbf{4.\ The geometric approach}}
\par\vspace{10pt}
 As usual we say that two functions $f,g$ (as well as their gaphs) have a \textsl{first-order contact} at a point $t_0$ if $f(t_0)=g(t_0)$ and $f'(t_0)=g'(t_0)$ provided that $f, g$ are defined on a neighborhood of  $t_0$ and the involved derivatives exist as finite numbers. The following elementary fact will provide a basis for our discussion.
\par\vspace{5pt}\textbf{Lemma 4.1}\quad \textsl{Let} $\phi_1, \phi_2$ \textsl{be two functions such that }
$$
\phi_1, \phi_2\quad differentiable\;\, on\;\, an\;\, interval\;\, I\,,\leqno(4.1)
$$
$$
W(\phi_1(x), \phi_2(x))\neq 0\qquad \forall x\in I\,.\leqno(4.2)
$$
\par
\textsl{In particular $\phi_1, \phi_2$ may satisfy the basic asumptions} (3.1) \textsl{on the interval} $[T, x_{{}_0}[$. \textsl{If $f$ is differentiable on $I$ then for each $t_0\in I$ there exists a unique function in the family 
$\mathcal{F}:= span (\phi_1, \phi_2)$ having a first-order contact with $f$ at $t_0$. Denoting this function by $F^* (x; t_0)$ we have}
$$
F^*(x; t_0)= f_1^*(t_0) \phi_1(x)+f_2^*(t_0)\phi_2(x)\;\,,\quad x\in I\,,\leqno(4.3)
$$
{\sl where}
$$\begin{cases}
f_1^*(t_0):=W(f, \phi_2; t_0)/W(t_0)\equiv (f(t)/\phi_2(t))'/(\phi_1(t)/\phi_2(t))'\big|_{t=t_0}\ ,\cr\cr
f_2^*(t_0):=-W(f, \phi_1; t_0)/W(t_0)\equiv (f(t)/\phi_1(t))'/(\phi_2(t)/\phi_1(t))'\big|_{t=t_0}\ .\end{cases}\leqno(4.4)$$
\textsl{If $f\in \mathcal{F}$ then $F^*(x; t_{{}_0})\equiv f(x)$ for any chosen to.}
\par\vspace{5pt}\textbf{Definition 4.1}\quad \textsl{In the quantity $F^*(x; t_0)$ we fix $x\in I$, say $x=T$, and consider the function}
$$
F^*(t):=F^*(T; t)\equiv \phi_1(T) f_1^*(t)+\phi_2(T) f^*_2(t)\,,\quad t\in I\,,\leqno(4.5)
$$
\textsl{which we call the contact indicatrix of order one of the function $f$ at the point $t$ with respect to the family $\mathcal{F}$ and the straight line} $x=T$.
\par
In the sequel we always suppose $I=[T, x_{{}_0}[$; the choice $x=T$ is merely a matter of convenience; any vertical line intersecting the $x$-interval $I$ can do the same. $F^*(t)$ represents the ordinate of the point of intersection between the vertical line $x=T$ and the curve $y=f^*_1(t)\phi_1(x)+f_2^*(t)\phi_2(x)$ where $t$ is thought of as fixed. By (4.2) $\phi_1$ and $\phi_2$ do not vanish simultaneously hence $F^*$ is a nontrivial linear combination of $f^*_1, f^*_2$. It may happen that, for some choices of $T$, $F^*$ coincides with $f^*_1$ or $f^*_2$, a constant factor apart, according as $\phi_2(T)=0$ or $\phi_1(T)=0$; this simply means that in a particular situation $F^*$ may be a redundant quantity, otherwise $F^*$ has its own pregnant geometric meaning. For instance if $(\phi_1, \phi_2)\equiv (x, 1)$ and $I=[0, +\infty)$ then $F^*(0; t)\equiv f^*_2(t)$. This cannot happen however if condition (3.1)$_3$, explicitly assumed as a matter of convenience, is satisfied.
\par
Using (4.4) $F^*$ may be represented as 
$$
F^*(x)=\displaystyle{1\over W(x)}[\phi_1(T)\cdot W(f, \phi_2; x)-\phi_2(T) W(f, \phi_1; x)]=\leqno(4.6)
$$
$$
=\displaystyle{1\over W(x)}\cdot W (f(x), \phi_1(T) \phi_2(x)-\phi_2(T) \phi_1(x))\equiv W(\Phi(x), f(x))/W(x)
$$
where we have put 
$$
\Phi(x):=\phi_2(T) \phi_1(x)-\phi_1(T) \phi_2(x)\,.\leqno(4.7)
$$
\par
If $f\in AC([T, x_{{}_0}[)$, then the function $F^*$ is defined almost everywhere on $I$ and is Lebesgue-summable on every compact interval of $I$.
\par\vspace{5pt}\textbf{Lemma 4.2}\quad (Representations of $f$ in terms of $F^*$, $f_1^*, f_2^*$). {\sl Under the basic assumptions} (3.1), {\sl except possibly} (3.1)$_2$, {\sl and with the foregoing notations let $f\in AC([[T, x_{{}_0}[)$. Then} 
$$
\Phi(x)\neq 0\quad\forall x \in ]T, x_0]\,;\leqno(4.8)
$$
$$
f(x)=c\Phi(x)+\Phi(x)\cdot \int\limits_{T_0}^x W(t) \Phi^{-2}(t) F^*(t) dt,\; x\in ]T, x_0[\,,\leqno(4.9)
$$
{\sl where $T_0$ is an arbitrarily fixed point in $]T, x_{{}_0}[$ and $c$ a suitable constant. Also, the following two representations are valid on the whole interval $[T, x_{{}_0}[$ with suitable constants}
$$
f(x)=c_2\phi_2(x)-\phi_2(x)\cdot\int\limits_T^x W(t) (\phi_2(t))^{-2} f^*_1(t) dt\,,\quad x\in [T, x_0[\,;\leqno(4.10)
$$
$$
f(x)=c_1\phi_1(x)+\phi_1(x)\cdot \int\limits_T^x W(t) (\phi_1(t))^{-2} f^*_2(t) dt\,,\quad x\in [T, x_0[\,.\leqno(4.11)
$$
\par\vspace{5pt}
\textbf{Proof.}\quad $\Phi$ is a nontrivial linear combination of $\phi_1, \phi_2$ as $\phi_i(T)\neq 0$ $(i=1, 2)$; hence, as mentioned after formula (3.3), $\Phi$ has at most one zero, namely $x=T$, and (4.8) follows. On $]T, x_{{}_0}[$ representation (4.6) can be written as
$$
F^*(x)={\Phi^2(x)\over W(x)}\left({f(x)\over \Phi(x)}\right)'\qquad a.e.\, on\,\; ]T, x_0]\,,\leqno(4.12)
$$
from whence (4.9) follows as $W\cdot\Phi^{-2}$ is continuous and $f/\Phi$ is absolutely continuous. Representations (4.10), (4.11) are similarly obtained and are valid on $[T, x_{{}_0}[$ by (3.1)$_3$.
\par
\hfill{$\Box$}\par\vspace{5pt}
So far the growth-order relation (3.1)$_2$ has played no role but it will play an essential one in obtaining our main results.
\par
We shall characterize the contingencies wherein each of the functions $f^*_i(x)$ and $F^*(x)$ admits of a finite limit as $x\to x_{{}_0}$ by means of suitable pairs of asymptotic expansions of $f$ and $f'$ with respect to the asymptotic scale $(\phi_1, \phi_2)$.
\par\vspace{5pt}
\textbf{Hypotheses for the three theorems in this section:}
$$\begin{cases}
\textnormal{(i)}& \textsl{the basic assumptions} (3.1) \textsl{about the pair of comparison}\cr
                & \textsl{functions} \phi_1, \phi_2\,;\cr
\textnormal{(ii)}& \textsl{a function} f\in AC([T, x_{{}_0}[)\,;\cr
\textnormal{(iii)}& \textsl{notations as in lemmas} 4.1-4.2 \textsl{and definition} 4.1\,.\cr\end{cases}\leqno(4.13)$$
\par\vspace{5pt}
\textbf{Theorem 4.3}\quad (The contingency: $\lim_{x\to x_{{}_0}} f^*_1(x)=a_1$). (I) \textsl{The following are equivalent properties:}
\par
(i) \textsl{There exists a finite limit}
$$
\lim_{x\to x_0} f^*_1(x)\equiv a_1\,.\leqno(4.14)
$$
\par
(ii) \textsl{It holds true the asymptotic relation}
$$
\left({f(x)\over\phi_2(x)}\right)'=a_1\left({\phi_1(x)\over\phi_2(x)}\right)'+o\left(\left({\phi_1(x)\over\phi_2(x)}\right)'\right)\,,\quad x\to x_0\,.\leqno(4.15)
$$
\par (iii)
\textsl{It holds true the pair of asymptotic relations}
$$\begin{cases}
f(x)=a_1\phi_1(x)+o(\phi_1(x))\cr\hspace{220pt},x\to x_0\,.\cr
\left(\displaystyle{f(x)\over\phi_2(x)}\right)'=a_1\left(\displaystyle{\phi_1(x)\over
\phi_2(x)}\right)'+o\left(\left(\displaystyle{\phi_1(x)\over\phi_2(x)}\right)'\right) \end{cases}\leqno(4.16)
$$
\par
\textsl{The above equivalences are simple consequences of the identity}
$$
(f(x)/\phi_2(x))'=(\phi_1(x)/\phi_2(x))' f_1^*(x)\,,\quad x\in [T, x_0[\,.\leqno(4.17)
$$
\par
\textsl{The constant} $a_1$ \textsl{in} (4.14), (4.15), (4.16) \textsl{is the same.}
\par
(II)\quad \textsl{The pair of conditions}
$$
\lim_{x\to x_0} f^*_1(x)=a_1\,;\quad \int\limits_T^{x_0}
\left(\displaystyle{\phi_1(t)\over\phi_2(t)}\right)'\,[f^*_1(t)-a_1] dt \;convergent\leqno(4.18)
$$
\textsl{is equivalent to the pair of asymptotic relations}
$$\begin{cases}
f(x)=a_1\phi_1(x)+a_2\phi_2(x)+o(\phi_2(x))\cr\hspace{220pt}, x\to x_0.\cr
\left(\displaystyle{f(x)\over\phi_2(x)}\right)'=a_1\left(\displaystyle{\phi_1(x)\over
\phi_2(x)}\right)'+o\left(\left(\displaystyle{\phi_1(x)\over\phi_2(x)}\right)'\right)\end{cases}\leqno(4.19)$$
\par
\textsl{The constant} $a_1$ \textsl{in } (4.18), (4.19) \textsl{is the same whereas } $a_2$ \textsl{is another suitable constant.}
\par
\textsl{If this is the case we have representation}
$$
f(x)=a_1\phi_1(x)+a_2\phi_2(x)-\phi_2(x)\cdot \int\limits_x^{x_0}
\left(\displaystyle{\phi_1(t)\over\phi_2(t)}\right)'\,[f^*_1(t)-a_1] dt, x\in [T, x_0[\,.\leqno(4.20)
$$
\par
The following intermediary result is an essential step in proving the subsequent main theorem.
\par\vspace{5pt}
\textbf{Theorem 4.4}\quad (Characterizations of a two-term asymptotic expansion). \textsl{Under assumptions} (4.13) \textsl{the following are equivalent properties:}
\par
(i) \textsl{It holds true an asymptotic expansion}
$$
f(x)=a_1\phi_1(x)+a_2\phi_2(x)+o(\phi_2(x))\,,\quad x\to x_0\,.\leqno(4.21)
$$
\par
(ii) \textsl{There exists a finite limit}
$$
\lim_{x\to x_0} \; \displaystyle{\phi_1(x)\over\phi_2(x)}\cdot \displaystyle{\int\limits_x^{x_0}}
W(t) (\phi_1(t))^{-2} f^*_2(t) dt
\equiv \displaystyle{\lim_{x\to x_0}}\; \displaystyle{\phi_1(x)\over\phi_2(x)}\cdot \displaystyle{\int\limits_x^{x_0}}\left({\phi_2(t)\over\phi_1(t)}\right)'f^*_2(t) dt\equiv -m\,.\leqno(4.22) $$

(iii) \textsl{There exists a finite limit}
$$
\lim_{x\to x_0} \displaystyle{\Phi(x)\over\phi_2(x)}\cdot \int\limits_x^{x_0}W(t) (\Phi(t))^{-2} F^*(t) dt\equiv -\displaystyle{l\over\phi_2(T)}\,.\leqno(4.23)
$$
\par
\textsl{If this is the case we have the following two representations}
$$\begin{cases} f(x)=a_1\phi_1(x)+a_2\phi_2(x)-\phi_1(x)\cdot\displaystyle{\int\limits_x^{x_0}W(t) (\phi_1(t))^{-2}[f^*_2(t)-m]} dt\equiv\\ \\
\equiv a_1\phi_1(x)+a_2\phi_2(x)-\phi_1(x)\cdot\displaystyle{\int\limits_x^{x_0}
\left(\displaystyle{\phi_2(t)\over\phi_1(t)}\right)'}[f^*_2(t)-m] dt\,,\quad x\in [T, x_0[\,;\end{cases}\leqno(4.24)$$
$$
f(x)=a_1\phi_1(x)+a_2\phi_2(x)-\Phi(x)\cdot\int\limits_x^{x_0}W(t) (\Phi(t))^{-2}[F^*(t)-l] dt\,,\, x\in ]T, x_0[\,.\leqno (4.25)
$$
\par
The validity of (4.21) may be expressed by the geometric locution: ``\textsl{the graph of} $f$ \textsl{admits of the curve } $y=a_1\phi_1(x)+a_2\phi_2(x)$ \textsl{as an asymptotic curve in the family} $\mathcal{F}\equiv span (\phi_1, \phi_2)$, \textsl{as} $x\to x_{{}_0}$." At the end of this section we suggest an expressive way of reading theorem 4.4 which for the time being looks like a technical lemma.
\par
In general there is no immediate relationship between the numbers $a_i, m$ and $l$ as in the case discussed in the following theorem, one of the main results in the paper.
\par\vspace{5pt}
\textbf{Theorem 4.5} \quad(The contingency $\lim_{x\to x_{{}_0}} f^*_2(x)=a_2$; characterizations of a limit tangent curve).\textsl{Let assumptions} (4.13) \textsl{hold true.}
\par
(I)\quad \textsl{The following are equivalent properties }
\par
(i) \textsl{There exists a finite limit}
$$
\lim_{x\to x_0} F^*(x)\equiv\gamma\,.\leqno(4.26)
$$
\par
(ii) \textsl{There exists a finite limit}
$$
\lim_{x\to x_0} f^*_2(x)\equiv a_2\,\quad {\rm (see} \ (4.29)_2 \ {\rm below)}.\leqno(4.27)$$
\par
(iii) \textsl{The following two limits exist as finite numbers}
$$
\lim_{x\to x_0} f^*_1(x)\equiv a_1\quad\,;\quad\lim_{x\to x_0} f^*_2(x)\equiv a_2\,.\leqno(4.28)
$$
\par
(iv) \textsl{It holds the pair of asymptotic relations}
$$
f(x)=a_1\phi_1(x)+a_2\phi_2(x)+o(\phi_2(x))\,, x\to x_0\,,\leqno (4.29)_1
$$
$$
\left(\displaystyle{f(x)\over\phi_1(x)}\right)'=a_2\left(\displaystyle{\phi_2(x)\over\phi_1(x)}\right)'+o\left[\left(\displaystyle{\phi_2(x)\over\phi_1(x)}\right)'\right]\,,\quad x\to x_0\,.\leqno(4.29)_2
$$
\par
(v) \textsl{It holds the pair of asymptotic relations}
$$
f(x)=a_1\phi_1(x)+a_2\phi_2(x)+o(\phi_2(x))\,, x\to x_0\,,\leqno (4.30)_1
$$
$$
\left(\displaystyle{f(x)\over\phi_2(x)}\right)'=a_1\left(\displaystyle{\phi_1(x)\over\phi_2(x)}\right)'+o\left[\displaystyle{\phi_2(x)\over\phi_1(x)}\left(\displaystyle{\phi_1(x)\over\phi_2(x)}\right)'\right]\,,\quad x\to x_0\,.\leqno(4.30)_2
$$
\par
(vi) \textsl{There exists a function} $\overline F$, \textsl{Lebesgue-summable on every compact interval of} $I$ \textsl{such that}
$$
\overline F(x)=o(1)\,,\quad x\to x_0\,,\leqno(4.31)_1
$$
\textsl{and}
$$
f(x)=a_1\phi_1(x)+a_2\phi_2(x)-\Phi(x)\int\limits_x^{x_0}W(t) (\Phi(t))^{-2}\overline F(t) dt\,,\quad x\in ]T, x_0[\,.\leqno (4.31)_2
$$
\par
\textsl{If this is in the case then}
$$
\overline F(x)=F^*(x)-\gamma\quad a.e.\, on \, I\,.\leqno(4.32)
$$
\par
(vii) \textsl{There exists a function $\overline f_2$, Lebesgue-summable on every compact interval of $I$ such that}
$$
f_2(x)=o(1)\,,\quad x\to x_0\,,\leqno(4.33)_1
$$
\textsl{and}
$$
f(x)=a_1\phi_1(x) +a_2\phi_2(x)-\phi_1(x)\cdot \int_x^{x_0}\left({\phi_2(t)\over \phi_1(t)}\right)' \overline f_2(t)dt\,,\quad x\in [T, x_0[\,.\leqno(4.33)_2
$$
\par
\textsl{If this is the case then}
$$
\overline f_2(x)=f_2^*(x)-a_2\qquad a.e.\; on\; I\,.\leqno(4.34)
$$\par
(II)\quad \textsl{Whenever properties in part} (I) \textsl{hold true then:}
\par
(viii) \textsl{The family of curves whose equations with respect to cartesian coordinates} $x, y$ \textsl{are}
$$
y=F^*(x; \xi)\equiv f^*_1(\xi)\phi_1(x)+f^*_2(\xi)\phi_2(x)\,,\leqno(4.35)
$$
\textsl{admits of a} ``\textsl{limit position}" \textsl{as} $\xi\to x_{{}_0}$, \textsl{namely}
$$
y=a_1\phi_1(x)+a_2\phi_2(x)\,,\leqno(4.36)
$$
\textsl{whose right-hand side is an asymptotic expansion of $f$, as $x\to x_{{}_0}$, formally differentiable once in the sense of relation} (4.29)$_2$. \textsl{ We say that the graph of $f$ admits of a } ``\textsl{limit tangent curve}" \textsl{in the family} $\mathcal{F}$ \textsl{as $x\to x_{{}_0}$, and this is a stronger contingency than the existence of an asymptotic curve as in theorem} 4.4.
\par
(ix) \textsl{The numbers} $\gamma, a_1, a_2$ \textsl{are linked by relation }
$$
\gamma=a_1\phi_1(T)+a_2\phi_2(T)\,\leqno(4.37)
$$
\textsl{whereas for the numbers $m, l$ appearing in theorem } 4.4 \textsl{we have } $m=a_2, l=\gamma$.
\par
(x) \textsl{They hold the two representations} (4.24) \textsl{and }(4.25) \textsl{with} $m=a_2$ \textsl{and} $l=\gamma$.
\par\vspace{5pt}\textsl{Remarks.}\quad By (4.4) the limit relations (4.28), when written out explicitly, coincide respectively with the asymptotic relations (4.15) and (4.29)$_2$ which can be also written in equivalent forms as
$$
W(f, \phi_2; x)=a_1 W(x)+o(W(x))\,,\quad x\to x_0\,,\leqno(4.38)
$$
$$
W(f, \phi_1; x)=-a_2 W(x)+o(W(x))\,,\quad x\to  x_0\,.\leqno(4.39)
$$\par
By looking at the mere formal aspect it is not self-evident that (4.29)$_2$ is stronger than (4.15): this follows instead from the detailed results in the foregoing theorems.
\par
That the sole relation (4.29)$_2\equiv$ (4.27), implies (4.29)$_1$ is a trivial consequence of (3.2) whereas relation (4.30)$_2$, which is a reinforced form of (4.15), does  not generally imply (4.30)$_1$. A trivial counterexample is provided by 
$$
\phi_1(x):=x\,;\quad \phi_2(x):=1\,;\quad f(x):= x+\log (\log x)\,;\quad x_0=+\infty\,.
$$
\par
\centerline{*\quad *\quad *}
\par
Before closing this section we mention how theorem 4.4 can be given a more expressive asymptotic meaning. Let us notice that a quantity such as
$$
\phi(x)\cdot\displaystyle{\int\limits_x^{x_0}}(1/\phi(t))' f(t) dt\leqno(4.40)
$$
is a kind of ``\textsl{weighted integral mean of} $f$" and that its limit as $x\to x_{{}_0}$ can be considered, the sign apart, as a ``\textsl{generalized limit of} $f(x)$ \textsl{as} $x\to x_{{}_0}$" for the simple reason that a trivial application of L'Hospital's rule yields
$$
\lim_{x\to x_0}\displaystyle{\displaystyle{\int\limits_x^{x_0}}(1/\phi(t))' f(t) dt\over (1/\phi(x))}=-\lim_{x\to x_0} f(x)
$$
provided that: $\phi\in AC([T, x_{{}_0}[)$; $\phi'(x)>0$ or $\phi'(x)<0$ a.e.; $\lim_{x\to x_{{}_0}}\phi(x)=\pm\infty$ and $\lim_{x\to x_{{}_0}} f(x)$ exists in $\overline{\mathbb R}$.
\par
In the very simple case $\phi(x)\equiv x$ on $[T,+\infty)$ (4.40) reduces (the sign apart) to
$$
\lim_{x\to +\infty} x\cdot\int\limits_x^{+\infty} t^{-2} f(t) dt\,.\leqno(4.41)
$$
\par
A result by Ostowski [12; IV, pp. 65-68] states that the limit (4.41) is equivalent to the simpler limit
$$
\lim_{x\to +\infty}  {1\over x} \cdot\int\limits_T^x f(t) dt\,.\leqno(4.42)
$$
\par
This limit appears here and there in the literature in problems related to the asymptotic behavior of solutions to ordinary differential equations, to the asymptotic behavior of Laplace transform and so on. In the case that $f$ is $p$-periodic the limit (4.42) exists and equals the usual mean of $f$ ``${1\over p}\int_T^{T+p} f$". In general if the quantity (4.42) is defined in $\mathbb{R}$ it may called the ``\textsl{asymptotic mean of} $f$ \textsl{at} $+\infty$". By analogy we may label the quantity
$$
-\lim_{x\to x_0}\phi(x)\cdot\int\limits_x^{x_0} (1/\phi(t))' f(t) dt\,,\leqno(4.43)
$$
if it is defined as a real number and with the above-specified restrictions on $\phi$, by the locution ``\textsl{asymptotic mean of} $f(x)$, \textsl{as} $x\to x_{{}_0}$, \textsl{with respect to the weight function} $\phi$". Of course some regularity condition on $f$ is required to give meaning to the foregoing integrals. With this terminology the equivalence ``(4.21)$\Leftrightarrow$(4.22)" may be reformulated as follows:
\par\vspace{5pt}
\textbf{Theorem 4.4 reformulated.}\quad \textsl{Under assumptions} (4.13) \textsl{the function $f$ admits of an asymptotic expansion with respect to the asymptotic scale} $(\phi_1, \phi_2)$, \textsl{as} $x\to x_{{}_0}$, \textsl{iff the associated geometric quantity} $f^*_2$ \textsl{admits of an asymptotic mean with respect to  the weight function} $\phi_1(x)/\phi_2(x)$ \textsl{as} $x\to x_{{}_0}$.
\par
By further investigation it is found out that the generalized asymptotic mean (4.43) is equivalent to the standard and simpler asymptotic mean (4.42) whenever $\phi$ is regularly varying at $+\infty$. We shall not go into the details of this subject in this paper.
\par\vspace{20pt}
\centerline{\textbf{5.\ The factorizational approach. Estimates of the remainder}}
\par\vspace{10pt}
In this section our basic assumptions on $\phi_1, \phi_2$ are the strenghtened ones, i.e. (3.1) plus (3.4), and $f\in AC^1(I)$. The treatment is based on canonical factorizations of the operator $L$ defined by (3.5)-(3.6). The link between the formulas in this section and those in the preceding one is provided by the following simple fact.
\par\vspace{5pt}
\textbf{Lemma 5.1.}\ {\it Hypotheses: {\rm (i)} the basic assumptions {\rm (3.1)} and {\rm (3.4);\ (ii)}  $f\in AC^1(I);$ {\rm (iii)} let $\phi(x):=c_1\phi_1(x)+c_2\phi_2(x)$ be such that $\phi(x)\neq 0$\ on\ some\ interval\ $J\subset I;$ {\rm (iv)} let}
$$
\widetilde f(x):=W(f(x), \phi(x))/W(x)\,.\leqno(5.1)
$$
\par
{\it Thesis}: \textsl{it holds the formula}
$$
\widetilde f'(x)=-\phi(x) (W(x))^{-1}\cdot L [f(x)]\quad a.e.\;\, on\;J\,,\leqno(5.2)
$$
\textsl{where $L$ is the differential operator} {\rm (3.5)-(3.6)}, \textsl{whence a representation of type}
$$
\widetilde f(x)=c-\int\limits_T^x \phi(t) (W(t))^{-1}\cdot L [f(t)]dt\,,\quad x\in J\,.\leqno(5.3)
$$
\par
\textsl{Replacing $\widetilde f$ by any of the three functions} $f^*_1, f^*_2, F^*$ \textsl{defined in} \S 4 \textsl{we get representations}
$$
f^*_1(x)=\overline c_1-\int\limits_T^x \phi_2(t) (W(t))^{-1}\cdot L[f(t)]dt\,,\quad x\in [T, x_0[\,;\leqno(5.4)
$$
$$
f^*_2(x)=\overline c_2+\int\limits_T^x \phi_1(t) (W(t))^{-1}\cdot L[f(t)]dt\,,\quad x\in [T, x_0[\,;\leqno(5.5)
$$
$$
 F^*(x)=\overline c+\int\limits_T^x \Phi(t) (W(t))^{-1}\cdot L[f(t)]dt\,,\quad x\in [T, x_0[\,.\leqno(5.6)
$$
\par\vspace{5pt}
\textbf{Proof.}\ Factorization (3.10) may be rewritten as 
$$
L[u]\equiv {W(x)\over\phi(x)}\left[{W(\phi(x), u)\over W(x)}\right]'\leqno(5.7)
$$
from whence (5.2) follows.
\par\hfill{$\Box$}
\par\vspace{5pt}
\textbf{The approach based on the canonical factorization of type (I) at $x_{{}_0}$.} The ``unique'' factorization of type (I) at $x_{{}_0}$ is
$$
L[u]\equiv {W(x)\over\phi_2(x)}\left[{(\phi_2(x))^2\over W(x)}\left({u\over\phi_2(x)}\right)'\right]'\,,\leqno(5.8)
$$
which gives rise to representation
$$
f(x)=c_1\phi_1(x)+c_2\phi_2(x)+\leqno(5.9)
$$
$$
+\phi_2(x)\int\limits_T^x (\phi_2(t))^{-2} W(t) dt\int\limits_T^t \phi_2(s)(W(s))^{-1}\cdot L[f(s)]ds\,,\quad x\in [T, x_0[\,.
$$
\par\vspace{5pt}
\textbf{The approach based on a canonical factorization of type (II) at $x_{{}_0}$.} By lemma 3.2 the simplest choice of a factorization of type (II) at $x_{{}_0}$, in terms of the given function $\phi_1, \phi_2$ is
$$
L[u]\equiv {W(x)\over\phi_1(x)}\left[{(\phi_1(x))^2\over W(x)}\left({u\over\phi_1(x)}\right)'\right]'\,,\leqno(5.10)
$$
which gives rise to representation
$$
f(x)=c_1\phi_1(x)+c_2\phi_2(x)+\leqno(5.11)
$$
$$
+\phi_1(x)\int\limits_T^x (\phi_1(t))^{-2} W(t) dt\int\limits_T^t \phi_1(s)(W(s))^{-1}\cdot L[f(s)]ds\,,\quad x\in [T, x_0[\,.
$$
\par
Representations (5.4)-(5.5) give the geometric meanings of the inner integrals appearing respectively in (5.9) and (5.11); if these inner integrals are replaced by (5.4) and (5.5) we get representations practically equivalent to (4.10) and (4.11).
\par
Representation (5.9), i.e. (4.10), is convenient to characterizing asymptotic expansions for $f$, matched to an asymptotic relation involving $(f/\phi_2)'$, see theorem 4.3, whereas representation (5.11), i.e. (4.11), is better fit to studying expansions for $f$ matched to relations involving $(f/\phi_1)'$ as in theorems 4.4, 4.5. The factorizational approach yields three integral representations which allow easy characterizations of certain asymptotic expansions through integral conditions involving $L[f]$, and this is the practical usefulness of this approach.
\par\vspace{5pt}
\textbf{Theorem 5.2}\quad (Restatement of conditions appearing in theorems 4.3-4.5). {\sl Assumptions are}: (3.1) \textsl{plus} (3.4) \textsl{and } $f\in AC^1(I)$.
\par
(I) (Refer to theorem 4.3). \textsl{Condition} (4.14) \textsl{for some real number $a_1$ is equivalent to the integral condition}
$$
\int\limits_T^{x_0} \phi_2(t)(W(t))^{-1}\cdot L[f(t)]dt\quad convergent\,,\leqno(5.12)
$$
\textsl{from whence it follows representation}
$$
f^*_1(x)=a_1+\int\limits_x^{x_0} \phi_2(t)(W(t))^{-1}\cdot L[f(t)]dt\,,\quad x\in [T, x_0[\,.\leqno(5.13)
$$
\par
\textsl{The pair of conditions } (4.18) \textsl{holds true iff the iterated improper integral}
$$
\int\limits_T^{x_0}\left({\phi_1(t)\over\phi_2(t)}\right)' dt \int\limits_t^{x_0} \phi_2(s)(W(s))^{-1}\cdot L[f(s)]ds\quad converges\,.\leqno(5.14)
$$
\par
(II) (Refer to theorem 4.5). \textsl{Condition} (4.27) \textsl{for some real number $a_2$ is equivalent to condition}
$$
\int\limits_T^{x_0} \phi_1(t) (W(t))^{-1}\cdot L[f(t)]dt\quad convergent\,,\leqno(5.15)
$$
\textsl{from whence we get representation}
$$
f^*_2(x)=a_2-\int\limits_x^{x_0} \phi_1(t)(W(t))^{-1}\cdot L[f(t)]dt\,,\quad x\in [T, x_0[\,.\leqno(5.16)
$$
\par
\textsl{Condition} (4.26) \textsl{for some real number $\gamma$ is equivalent to}
$$
\int\limits_T^{x_0} \Phi(t)(W(t))^{-1}\cdot L[f(t)]dt\quad convergent\,,\leqno(5.17)
$$
\textsl{which yields representation}
$$
F^*(x)=\gamma -\int\limits_x^{x_0} \Phi(t)(W(t))^{-1}\cdot L[f(t)]dt\,,\quad x\in [T, x_0[\,.\leqno(5.18)
$$
\par
Representations (5.13), (5.16) and (5.18) may be substituted into (4.20), (4.24) and (4.25) respectively so obtaining numerically useful formulas reported in the theorem below.
\par\vspace{5pt}
\textbf{Theorem 5.3}\quad (Representations and estimates of the remainder). \textsl{Let} $f\in AC(I)$ \textsl{and let its graph admit of the curve} $y=a_1\phi_1(x)+a_2\phi_2(x)$ \textsl{as a limit tangent curve in the family} $\mathcal{F}$, \textsl{as} $x\to x_{{}_0}$. \textsl{Put }
$$
R(x):=f(x)-a_1\phi_1(x)-a_2\phi_2(x)\,.\leqno(5.19)
$$
\par
(I) Integral representations. \textsl{In the situation of theorem} (4.5) \textsl{we have}
$$
R(x)=-\phi_1(x)\cdot \displaystyle{\int\limits_x^{x_0}} W(t) (\phi_1(t))^{-2}[f_2^*(t)-a_2]dt\equiv\leqno(5.20)
$$
$$
\equiv-\phi_1(x)\cdot \displaystyle{\int\limits_x^{x_0}} (\phi_2(t)/\phi_1(t))' [f_2^*(t)-a_2]dt\,,\; x\in [T, x_0[\,;
$$
$$
R(x)=-\Phi(x)\cdot \int\limits_x^{x_0} W(t) (\Phi(t))^{-2}[F^*(t)-\gamma]dt\,,\quad x\in ]T, x_0]\,.\leqno(5.21)
$$
\par
\textsl{If $\phi_1, \phi_2, f$ satisfy the stronger assumptions in theorem} 5.2 \textsl{then}
$$
R(x)=\phi_1(x)\cdot \int\limits_x^{x_0}(\phi_2(t)/\phi_1(t))' dt \int\limits_t^{x_0}\phi_1 (s)
(W(s))^{-1}\cdot L[f(s)] ds\,,\; x\in [T, x_0[\,;\leqno(5.22)
$$
$$
R(x)=-\phi_2(x)\cdot \int\limits_x^{x_0}(\phi_1(t)/\phi_2(t))' dt \int\limits_t^{x_0}\phi_2 (s)
(W(s))^{-1}\cdot L[f(s)] ds\,,\; x\in [T, x_0[\,;\leqno(5.23)
$$
$$
R(x)=\Phi(x)\cdot\int\limits_x^{x_0}W(t)(\Phi(t))^{-2} dt \int\limits_t^{x_0}\Phi(s)
(W(s))^{-1}\cdot L[f(s)] ds\,,\; x\in ]T, x_0[\,.\leqno(5.24)
$$
\par
(II) Estimates. \textsl{From} (5.20) \textsl{and } (5.21) \textsl{we get respectively}
$$
\vert R(x)\vert\leq \vert\phi_2(x)\vert\cdot \left( 
\mathrel{\mathop{ess.sup.}_{x<t<x_0}}
 \vert f_2^*(t)-a_2\vert\right)\,,\quad x\in [T, x_0[\,;\leqno(5.25)
$$
$$
\vert R(x)\vert\leq \left\vert{\phi_2(x)\over \phi_2(T)}\right\vert\cdot \left( 
\mathrel{\mathop{ess.sup.}_{x<t<x_0}}
 \vert F^*(t)-\gamma\vert\right)\,,\quad x\in [T, x_0[\,;\leqno(5.26)
$$
\par
\textsl{whereas from} (5.22) \textsl{and } (5.24) \textsl{we get respectively}
$$
\vert R(x)\vert\leq \vert\phi_2(x)\vert\cdot \int\limits_x^{x_0}\vert\phi_1(t)(W(t))^{-1}\cdot L[f(t)]\vert dt\,,\; x\in [T, x_0[\,;\leqno(5.27)
$$
$$
\vert R(x)\vert\leq \left\vert{\phi_2(x)\over \phi_2(T)}\right\vert\cdot \int\limits_x^{x_0}\vert\Phi(t)(W(t))^{-1}\cdot L[f(t)]\vert dt\,,\; x\in ]T, x_0[\,.\leqno(5.28)
$$
(III) Lagrange-type representations. \textsl{If $f\in C^1(I)$ then for each $x\geq T$ there exists} $\xi_1\in\overline{\mathbb R}$, $x\leq\xi_1\leq x_{{}_0}$, \textsl{such that}
$$
R(x)=\phi_2(x) (f^*_2(\xi_1)-a_2)\;\,,\quad (f_2^*(x_0):=a_2)\,;\leqno(5.29)
$$
\textsl{and for each} $x > T$ \textsl{there exists} $\xi_2\in\overline{\mathbb{R}}, x\leq\xi_2\leq x_{{}_0}$, \textsl{such that}
$$
R(x)={\phi_2(x)\over \phi_2(T)} (F^*(\xi_2)-\gamma)\;\,,\quad (F^*(x_0):=\gamma)\,.\leqno(5.30)
$$
\par
\textsl{Under the stronger assumptions in theorem} 5.2, (5.29) \textsl{and} (5.30) \textsl{may be respectively written as}
$$
R(x)=-\phi_2(x) \cdot \int\limits_{\xi_2}^{x_0}\phi_1(t)(W(t))^{-1}\cdot L[f(t)] dt\,,\; x\in ]T, x_0[\,;\leqno(5.31)
$$
$$
R(x)=-{\phi_2(x)\over\phi_2(T)} \cdot \int\limits_{\xi_2}^{x_0}\Phi(t)(W(t))^{-1}\cdot L[f(t)] dt\,,\; x\in ]T, x_0[\,,\leqno(5.32)
$$
\textsl{with the obvious agreement} $\displaystyle{\int_{x_{{}_0}}^{x_{{}_0}}}=0$.
\par\vspace{5pt}\textsl{Remark.}\quad Representation (5.23) comes out from (4.20) which holds true under conditions (4.18) and the weaker assumptions in theorem 4.3 granting that the curve $y=a_1\phi_1(x)+a_2\phi_2(x)$ is an asymptotic curve for the graph of $f$, as $x\to x_{{}_0}$, but not necessarily a limit tangent curve. However numerical estimates obtained from (5.23) cannot have the simple forms reported above due to the divergence of the integral $\int^{x_{{}_0}}(\phi_1/\phi_2)'$.
\par\vspace{20pt}
\centerline{\bf 6. The case of generalized convex functions}
\par\vspace{10pt}
The main result in this section states that: \textsl{If $f$ is a generalized convex function with respect to the Chebyshev system $(\phi_1, \phi_2)$ then the existence of an asymptotic expansion} (4.21) \textsl{automatically implies the existence of a limit tangent curve in the family} $\mathcal F$ \textsl{as} $x\to x_{{}_0}$ (theorem 4.5) \textsl{and this last contingency is even implied by the weaker relation}
$$
f(x)=a_1\phi_1(x)+O(\phi_2(x))\,,\; x\to x_0\,.
$$
\par
Before stating the precise result we point out that the asymptotic properties discussed so far as well as those to be discussed in this section do not depend in themselves on the signs of $\phi_1, \phi_2, W(\phi_1, \phi_2)$, whereas a decisive role in the concept of convexity is played by the monotonicity of certain functions and the types of monotonicity do depend on certain signs. So it is better to make a definite agreement about the signs and this is contained in the following standard definitions of Chebyshev systems and generalized convex functions.
\par\vspace{5pt}\textbf{Definition 6.1}\quad (Two-dimensional Chebyshev systems). \textsl{Let $(\psi_1, \psi_2)$ be an ordered pair of continuous functions on an interval} $J\in\mathbb{R}$.
\par
(I). $(\psi_1, \psi_2)$ \textsl{is a $T$-system} ($\equiv$ \textsl{Chebyshev system}) \textsl{on $J$ iff}
$$
\left\vert\begin{array}{ccc}
\psi_1(t_1)&\psi_1(t_2)\\
\psi_2(t_1)&\psi_2(t_2)
\end{array}
\right\vert >0\quad \forall t_1, t_2\in J\,\,;\quad t_1<t_2\,.\leqno(6.1)
$$
\par
(II). $(\psi_1, \psi_2)$ \textsl{is a $CT$-system} ($\equiv$ \textsl{complete Chebyshev system}) \textsl{on $J$ iff, in addition to} (6.1),
$$
\psi_1(t)>0\qquad \forall t\in J\,,\leqno(6.2)
$$
\textsl{without any a priori restriction on the sign of } $\psi_2$.
\par
(III). $(\psi_1, \psi_2)$ \textsl{is an $ET$-system} ($\equiv$ \textsl{extended Chebyshev system}) \textsl{on $J$ iff, in addition to} (6.1),
\textsl{the following two conditions are satisfied}:
$$
\psi_1, \psi_2\in C^1(J)\,;\leqno(6.3)
$$
$$
\left\vert\begin{array}{ccc}
\psi_1(t)&\psi'_1(t)\\ \\
\psi_2(t)&\psi'_2(t)
\end{array}
\right\vert >0\quad \forall t\in J\,.\leqno(6.4)
$$
\par
(IV). $(\psi_1, \psi_2)$ \textsl{is an $ECT$-system} ($\equiv$ \textsl{extended  complete Chebyshev system}) \textsl{on $J$ iff all conditions } (6.1), (6.2), (6.3), (6.4) \textsl{are satisfied}.
\par
The above locutions are those in the book by Karlin and Studden [10; chp. I]. It is known that $(\psi_1, \psi_2)$ is a $T$-system on $J$, except possibly for the sign of $\psi_2$, iff any nontrivial linear combination of $\psi_1, \psi_2$ has at most one zero on $J$. This is stated for general $T$-systems and for $J$ a compact interval in [10; th. 4.1, p.22], but a rereading of the proof shows that the argument remains unchanged for any interval. It is also known that in definition 6.1-(IV) condition (6.1) is redundant in so far as the three conditions (6.2), (6.3), (6.4) imply (6.1): see Coppel [1; prop. 5, p.89] or Karlin and Studden [10; th. 1.1, p. 376].
\par\vspace{5pt}\textbf{Definition 6.2}\quad \textsl{Let $(\psi_1, \psi_2)$ be a Chebyshev system on an interval $J$; a function $f: J\to\mathbb{R}$ is termed} ``\textsl{convex}" \textsl{on $J$ with respect to the system} $(\psi_1, \psi_2)$ \textsl{iff}
$$
U\left(\begin{array}{ccc}
\psi_1,&\psi_2,&f\\ \\
t_1,&t_2,&t_3
\end{array}
\right) :=\left\vert\begin{array}{ccc}
\psi_1(t_1)&\psi_1(t_2)&\psi_1(t_3)\\ \\
\psi_2(t_1)&\psi_2(t_2)&\psi_2(t_3)\\ \\
f(t_1)&f(t_2)&f(t_3)
\end{array}\right\vert \geq  0\leqno(6.5)
$$
\textsl{for each choice of} $t_1, t_2, t_3 \in J: t_1< t_2<t_3$. \textsl{It is termed} ``\textsl{strictly convex}" \textsl{iff the strict sign prevails in} (6.5). \textsl{Whenever} (6.5) \textsl{is satisfied we use notation} $f\in \mathcal C (\psi_1, \psi_2; J)$.
\par
As a standard reference for this class of functions (with respect to an $n$-dimensional Chebyshev system) we again quote Karlin and Studden [10; chp. XI].
\par
From now on in this section we shall be considering a pair of comparison functions which, besides satisfying all assumptions (3.1), form an $ECT$-system. According to the remarks following definition 6.1 it is enough to consider a pair $(\phi_1, \phi_2)$ satisfying
$$
\phi_1, \phi_2\in C^1(I)\; \,;\; I :=[T, x_0[\,;\leqno(6.6)_1
$$
$$
\phi_1(x)>>\phi_2(x)\; \,,\; x\to x_0^-\,.\leqno(6.6)_2
$$
$$
\phi_1(x)>0\; \forall x\in I\;\,;\; \phi_2(x)\neq 0\; \forall x\in I\,;\leqno(6.6)_3
$$
$$
W(x)\equiv W(\phi_1(x), \phi_2(x))>0\quad\forall x\in I\,.\leqno(6.6)_4
$$
\par
Our main results are collected in the following two theorems: the first one dealing with monotonicity properties and the second one with asymptotic properties of two-dimensional generalized convex functions.
\par\vspace{5pt}
\textbf{Theorem 6.1}\quad (Monotonicity properties of two-dimensional generalized convex functions). \textsl{If the pair $(\phi_1, \phi_2)$ satisfies all conditions} (6.6) \textsl{and if} $f\in\mathcal C(\phi_1, \phi_2; ]T, x_{{}_0}[)$ \textsl{then all the following properties hold true}:
\par
(i) $f\in AC]T, x_{{}_0}[$.
\par
(ii) \textsl{The three functions $f^*_1, f^*_2, F^*$ defined by} (4.4)-(4.5) \textsl{are defined a.e. on $]T, x_{{}_0}[$ and are monotonic on $]T, x_{{}_0}[\;\setminus N$ where $N$ is some Lebesgue null-set}: $f^*_2$ \textsl{is increasing}; $f_1^*$ \textsl{and} $F^*$ \textsl{have opposite types of monotonicity and} $(sign \phi_2)\cdot F^*$ \textsl{is increasing}.
\par
(iii) \textsl{Any function of type $[f(x)+a_1\phi_1(x)+a_2\phi_2(x)]/\widetilde\phi(x)$ is either constant or strictly monotonic on a suitable deleted neighborhood of any of the endpoints $T, x_{{}_0}$. Here $a_k$ are any constants and $\widetilde \phi$ is any nontrivial linear combination of} $\phi_1, \phi_2$.
\par 
(iv) \textsl{The following two limits exists simultaneously in $\overline{\mathbb R}$ and are equal}
$$
\lim_{x\to x_0} f(x)/\phi_1(x)=\lim_{x\to x_0}{(f(x)/\phi_2))'\over (\phi_1(x)/\phi_2(x))'}\,,\leqno(6.7)
$$
\textsl{the roles of} $\phi_1, \phi_2$ \textsl{being not interchangeable: compare with theorem 4.3}.
\par
Notice that the contingency $f\in \mathcal C(\phi_1, \phi_2; ]T, x_{{}_0}[)$ can be \underline{characterized} by the appropriate type of monotonicity of any of the functions $f^*_1, f_2^*, F^*$: see lemma 7.6 for a precise statement.
\par\vspace{5pt}
\textbf{Theorem 6.2}\quad (Asymptotic expansions of two-dimensional generalized convex functions). \textsl{Under the same assumptions on $\phi_1, \phi_2, f$ as in the foregoing theorem the following facts hold true}:
\par
(i) \textsl{We have the inference }
$$
f(x)=O(\phi_1(x))\,,\, x\to x_0\,\Rightarrow \begin{cases}
f(x)=a_1\phi_1(x)+o(\phi_1(x)),\cr\cr
(f(x)/\phi_2(x))'=a_1(\phi_1(x)/\phi_2(x))'
+o(\phi_1(x)/\phi_2(x))',\end{cases}\leqno(6.8)
$$
\textsl{for some constant $a_1$}: see theorem 4.3.
\par
(ii) \textsl{To the equivalent properties} (i)-(vi) \textsl{listed in theorem} 4.5-(I) \textsl{each of the following may be added}:
$$
F^*(x)=O(1)\,,\; x\to x_0\,;\leqno(6.9)
$$
$$
f^*_2(x)=O(1)\,,\; x\to x_0\,;\leqno(6.10)
$$
$$
f(x)=a_1\phi_1(x)+O((\phi_2(x))\,,\;  x\to x_0\quad for\;\, some\;\, constant \,\; a_1\,.\leqno(6.11)
$$
(Compare with properties (i), (ii) and (iv) in theorem 5.3).
\par
(iii) \textsl{Whenever} (4.26)-(4.28) \textsl{are satisfied then the following inequalities hold true}
$$
(sign\;\phi_1) (f^*_1(x)-a_1)\geq 0\,,\leqno(6.12)
$$
$$
(sign\;\phi_2) (F^*(x)-\gamma)\leq 0\,,\leqno(6.13)
$$
$$
f^*_2(x)-a_2\leq 0\,,\leqno(6.14)
$$
\textsl{for each } $x\in ]T, x_{{}_0}[\setminus N$ \textsl{where $N$ is as in theorem} 6.1, \textsl{and}
$$
R(x):=f(x)-a_1\phi_1(x)-a_2\phi_2(x)\geq 0\quad \forall x\in ]T, x_0[\,.\leqno(6.15)
$$
\par
\textsl{Moreover if there exists a point} $\xi\in ]T, x_o[$ \textsl{such that} $R(\xi)=0$ \textsl{then} $R(x)=0$ $\forall x\in [\xi, x_{{}_0}[$.
\par\vspace{5pt}\textsl{Remarks.}\ 1. The import of theorem 6.2 is that condition ``$f/\phi_1$ bounded" implies a one-term asymptotic expansion whereas ``$(f-a_1\phi_1)/\phi_2$ bounded" implies a two-term asymptotic expansion; moreover each expansion is formally differentiable in a suitable sense (not the same in the two cases).
\par
2. For practical applications it is important to bear in mind that under the stronger regularity conditions (3.1)-(3.4) and $f\in AC^1 (I)$ we have
$$
f\in \mathcal C(\phi_1, \phi_2; (I)\; \Leftrightarrow \; L[f(x)]\geq 0\;\, a.e.\;\, on \; I\,\leqno(6.16)
$$
where $L$ is the operator (3.5)-(3.6): see lemma 7.6. For a funtion $f$ satisfying $L[f]\geq 0$ the whole asymptotic theory developed so far admits of simpler proofs.
\par\vspace{20pt}
\centerline{\bf 7. Proofs.}
\par\vspace{10pt}
\textbf{Lemma 7.1}\quad (Trivial Wronskian identities). \textsl{If} $\phi\in \mathcal F :=span(\phi_1, \phi_2)$ \textsl{then:}
\par
(I)\quad \textsl{For each} $i=1,2$
$$
W(\phi, \phi_i; x)= c_i W(x)\quad\,,\quad x\in I\,,\leqno(7.1)
$$
\textsl{where} $c_i$ \textsl{is a suitable constant, possibly} $c_i=0$.
\par
(II)\quad \textsl{If} $\phi$ \textsl{is such that}  $\phi(x)\neq 0$ $\forall x\in J$, $J$ \textsl{a subinterval of} $I$, \textsl{then for at least one value of} $i=1,2$ \textsl{we have the identity}
$$
W(x) \phi^{-2}(x)=\overline c_i(\phi_i(x)/\phi(x))'\quad\,,\quad x\in J\,,\leqno(7.2)
$$
\textsl{for a suitable constant} $\overline c_i\neq 0$.
\par\vspace{15pt}
\textbf{Proof.}\quad For $\phi=a_1\phi_1+a_2\phi_2$ we have
$$
W(a_1\phi_1+a_2\phi_2, \phi_i)=
{{{i=1}\atop{\nearrow}}\atop{\searrow\atop{\scriptstyle{i=2}}}}
\kern-.4cm{\begin{array}{ll}
&a_2 W(\phi_2, \phi_1)=-a_2 W(\phi_1, \phi_2)\\ \\
&a_1 W(\phi_1, \phi_2)\,.
\end{array}}
\leqno(7.3)
$$
\par
If $\phi\neq 0$ on $J$ then at least one $a_i$ is non-zero hence, by (7.3), at least one $c_i$ in (7.1) is non-zero. For each such value of $i$ we may write
$$
W(x)={1\over c_i} W(\phi(x), \phi_i(x))\equiv {1\over c_i} \phi^2(x) (\phi_i(x)/\phi(x))'\,,\quad x\in J\,.\leqno(7.4)$$
\hfill{$\Box$}

\vspace{5pt}
\textbf{Lemma 7.2}\quad (Some properties of $F^*$).
\par (I)\quad \textsl{Function} $F^*$ \textsl{defined by} (4.4)-(4.5) \textsl{remains unchanged if we add to} $f$ \textsl{a function of type} $c$ $\Phi(x)$ \textsl{where } $\Phi$ \textsl{is defined by} (4.7).
\par (II)\quad \textsl{Condition}
$$
F^*(x)\equiv\gamma=\; constant\quad \forall x\in I\leqno(7.5)
$$
\textsl{holds true iff } $f\in\mathcal F := span(\phi_1, \phi_2)$.
\par\vspace{15pt}
{\it Proof.}\ (I) follows from (4.6). If $f=c_1\phi_1+c_2\phi_2$ then (4.6) and (3.20) imply (7.5). Viceversa if (7.5) holds true representation (4.9) gives
$$
f(x)=c\Phi(x)+\gamma\Phi(x)\cdot \int\limits_{T_0}^x W(t)\Phi^{-2}(t) dt \stackrel{by (7.2)}{=}
 \overline c\Phi(x)+\overline\gamma\phi_i(x)\,,\; x\in I\,,\leqno(7.6)
$$
for some $i=1,2$ and suitable constants $\overline c, \overline\gamma$; hence $f\in \mathcal F$.
\hfill{$\Box$}
\par\vspace{15pt}
Lemmas 7.1 and 7.2 do not depend on the asymptotic relation (3.1)$_2$; on the contrary much of the subsequent results depend on all four assumptions (3.1).
\par\vspace{15pt}
\textbf{Lemma 7.3}\quad (Some properties of $\Phi$). \textsl{Under the basic assumptions} (3.1), \textsl{the following are true}:
$$
\Phi(x)\sim\phi_2(T) \phi_1(x)\,,\quad x\to x_0\,;\leqno(7.7)
$$
$$
W(x)\Phi^{-2}(x)\sim {1\over(\phi_2(T))^2} (\phi_2(x)/\phi_1(x))'\,,\quad x\to x_0\,;\leqno(7.8)
$$
$$
\int\limits_x^{x_0} W(t)\Phi^{-2}(t) dt=-{1\over\phi_2(T)} (\phi_2(x)/\Phi(x))\,,\quad x\in ]T, x_0[\,;\leqno(7.9)
$$
$$
\int\limits_x^{x_0} W(t)\Phi^{-2}(t) dt\sim -{1\over(\phi_2(T))^2} (\phi_2(x)/\phi_1(x))\,,\quad x\to x_0\,.\leqno(7.10)
$$
\par\vspace{15pt}
{\it Proof.}\  Relation (7.7) is a direct conseguence of (4.7), (3.1)$_2$ and (3.1)$_3$. Relation (7.8) follows from (7.7):
$$
W(x)\Phi^{-2}(x)\sim W(x)(\phi_2(T)\phi_1(x))^{-2}\equiv (\phi_2(T))^{-2} (\phi_2(x)/\phi_1(x))'\,.
$$
\par
Relation (11.8) implies, by (3.2), the convergence of the integral appearing in (7.9); formula (7.9) follows from
$$
W(x)\Phi^{-2}(x)\stackrel{by (4.7)}{=}
{1\over \phi_2(T)} W(\Phi, \phi_2; x) \Phi^{-2}(x)\equiv \leqno(7.11)
$$
$$
\equiv {1\over \phi_2(T)} (\phi_2(x)/\Phi(x))'\,,\qquad x\in ]T, x_0[\,.
$$
\hfill{$\Box$}
\par\vspace{5pt}
{\it Proof of theorem} 4.3. (I). Let (4.14) hold true. From (4.4)
$$
f_1^*(x)\equiv W(f, \phi_2; x) (W(x))^{-1}=-\left({f(x)\over \phi_2(x)}\right)' (\phi_2(x))^2 (W(x))^{-1}\,,\leqno(7.12)
$$
from whence (4.17) follows and from (4.17) the equivalence ``(i)$\Leftrightarrow$ (ii)" is at once inferred, together with a representation of type
$$
f(x)/\phi_2(x)=c+\int\limits_T^x (\phi_1(t)/\phi_2(t))' f_1^*(t) dt\,,\quad x\in I\,.\leqno(7.13)
$$
\par
Moreover if (4.15) holds true then, by the divergence of the integral $\int^{x_{{}_0}}(\phi_1/\phi_2)'$, we infer
$$
\begin{array}{ll}
f(x)/\phi_2(x)&=\overline c+a_1(\phi_1(x)/\phi_2(x))+o(\phi_1(x)/\phi_2(x))=\\ \\
&=a_1(\phi_1(x)/\phi_2(x))+o(\phi_1(x)/\phi_2(x))\,,\quad x\to x_0\,,
\end{array}
$$
that is to say the first relation (4.16), hence ``(ii)$\Leftrightarrow$ (iii)".
\par
(II). If the first condition (4.18) is satisfied we have (7.13) which we may rewrite as
$$
f(x)/\phi_2(x)=\overline c+a_1{\phi_1(x)\over \phi_2(x)}+\int\limits_T^x (\phi_1(t)/\phi_2(t))' [f_1^*(t)-a_1] dt\,.\leqno(7.14)
$$\par
If the second condition (4.18) is also satisfied then we rewrite (7.14) as
$$
f(x)/\phi_2(x)=a_1{\phi_1(x)\over \phi_2(x)}+a_2-\int\limits_x^{x_0} (\phi_1(t)/\phi_2(t))' [f_1^*(t)-a_1] dt\,,\leqno(7.15)
$$
which implies representation (4.20) and
the first relation (4.19). The second relation (4.19) holds true by part (I) of the theorem. Viceversa if (4.19) hold true then, by part (I), we have (4.14), (7.13) and (7.14). The first relation (4.19) at once implies the second condition (4.18).
\hfill{$\Box$}
\par\vspace{15pt}
\textbf{Proof of theorem 4.4.}\quad (i)$\Leftrightarrow$ (iii): if (4.21) is true then $f=a_1\phi_1+o(\phi_1)$; using this relation and (7.7) into representation (4.9) we get condition
$$
\int\limits^{x_0} W \Phi^{-2} F^*\quad convergent\,,\leqno(7.16)
$$
and (4.9) may be rewritten as
$$
f(x)=a_1\phi_1(x)+c\phi_2(x)-\Phi(x)\cdot
\int\limits^{x_0}_x W \Phi^{-2} F^*\quad x\in ]T, x_0[\,,\leqno(7.17)
$$
where $a_1$ is the same constant as in (4.21) and $c$ is a suitable constant. From (4.21) and (7.17) we get condition (4.23) where we have denoted the value of the limit by $-l/\phi_2(T)$ for reasons of convenience. Viceversa if (4.23) holds true we have
$$
\Phi(x)\cdot\int\limits^{x_0}_x W \Phi^{-2} F^*=-{l\over\phi_2(T)}\phi_2(x)+o(\phi_2(x))\,,\leqno(7.18)
$$
and from (7.17) we get (4.21). Now, using (7.9) we rewrite (7.17) in the form
$$
f(x)=a_1\phi_1(x)+\left(c+{l\over\phi_2(T)}\right)\phi_2(x)-
\Phi(x)\cdot\int\limits^{x_0}_x W(t) (\Phi(t))^{-2} [F^*(t)-l],\ x\in ]T, x_0[\ ,
\leqno(7.19)$$
where the last term is
$$
\Phi(x)\cdot\int\limits^{x_0}_x W \Phi^{-2} (F^*-l)=\Phi(x)\cdot \int\limits^{x_0}_x W \Phi^{-2} F^*+{l\over\phi_2(T)}\phi_2(x)
\stackrel{by (7.18)}{=}\; o(\phi_2(x))\,.
$$
\par
Hence we get representation (4.25) with $a_2=c+{l\over \phi_2(T)}$.
\par
In a similar way we show the equivalence `` (i) $\Leftrightarrow$ (ii)" and (4.24) using representation (4.11) instead of (4.9).
\hfill{$\Box$}
\par\vspace{15pt}
{\it Proof of theorem} 4.5.\quad Part (I). First we notice that ``(ii)$\Leftrightarrow$ (iv)" as (4.29)$_2$ is an equivalent way of writing (4.27), by the very definition (4.4) of $f^*_2$, and (4.29)$_1$ is an automatic consequence of (4.29)$_2$ due to (3.2). The plan of our proof consists in proving the following inferences:
\par
(iv) $\Rightarrow$ (i) $\Rightarrow$ (vi) $\Rightarrow$ (iv); (ii) $\Leftrightarrow$ (iii); (vi) $\Rightarrow$ (v) $\Rightarrow$ (i).
\par
The equivalence ``(i) $\Leftrightarrow$ (vii)'' is perfectly analogous to ``(i) $\Leftrightarrow$ (vi)''.
\par
(iv) $\Rightarrow$ (i). As shown at the outset of the proof
of theorem 4.4 we have at our disposal representation (7.17) from which we get
$$
(f(x)/\phi_1(x))'\kern-2pt=\kern-2pt c(\phi_2(x)/\phi_1(x))'-\left({\Phi(x)\over \phi_1(x)}\right)'\cdot \kern-4pt\int\limits_x^{x_0} W\Phi^{-2} F^*
\kern-2pt+\kern-2pt{W(x)\over \phi_1(x)\Phi(x)} F^*(x)\kern-3pt=\leqno(7.20)
$$
$$
\stackrel{by (4.7)}{=}c(\phi_2(x)/\phi_1(x))'+\phi_1(T)(\phi_2(x)/\phi_1(x))'\cdot\int\limits_x^{x_0} W \Phi^{-2} F^*+
{W(x)\over \phi_1(x)\Phi(x)} F^*(x)=
$$
$$
\stackrel{by (7.7)}{=}c(\phi_2(x)/\phi_1(x))'+o((\phi_2(x)/\phi_1(x))')+
{1\over \phi_2(T)} (\phi_2(x)/\phi_1(x))' F^* (x) [1+o(1)]\,.
$$
\par
This trivially implies (4.26) by assumption (4.29)$_2$.
\par
(i) $\Rightarrow$ (vi). Relations (4.26) and (7.8) imply (7.16) and we may rewrite representation (4.9) in the form
$$
f(x)=\overline c\Phi(x)-\Phi(x)\cdot\int\limits_x^{x_0} W\Phi^{-2} F^*\equiv\leqno(7.21)
$$
$$
\equiv \overline c \Phi(x)-\gamma\Phi(x)\cdot \int\limits_x^{x_0} W\Phi^{-2}-\Phi(x)\cdot \int\limits_x^{x_0} W(t)\Phi^{-2}(t)[F^*(t)-\gamma]dt=
$$
$$
\stackrel{by (7.9)}{=}a_1\phi_1(x)+a_2\phi_2(x)-\Phi(x)\cdot \int\limits_x^{x_0} W(t)\Phi^{-2}(t) [F^*(t)-\gamma]dt\,,
$$
where $a_1, a_2$ are suitable constants. This is (4.37) which implies the assertion in (vi).
\par 
(vi) $\Rightarrow$ (iv).  From (4.31) and (7.9) we infer that
$$
\Phi(x)\cdot \int\limits_x^{x_0} W\Phi^{-2}\overline F=o(\phi_2(x))\,,\quad x\to x_0\,;\leqno(7.22)
$$
hence representation (4.31)$_2$ implies (4.29)$_1$ with the same constants $a_1, a_2$. Differentiating (4.31)$_2$ we get
$$
(f(x)/\phi_1(x))'=a_2(\phi_2(x)/\phi_1(x))'-\left({\Phi(x)\over\phi_1(x)}\right)'\cdot \int\limits_x^{x_0} W\Phi^{-2} \overline F
+{W(x)\overline F(x)\over \phi_1(x)\Phi(x)}=\dots =(4.29)_2\leqno(7.23)$$
by the same calculations in (7.20) with $F^*$ replaced by $\overline F$.
\par
(iii) $\Rightarrow$ (ii): obvious. (ii) $\Rightarrow$ (iii). As noticed at the outset of the present proof relation (4.27) is nothing but relation (4.29)$_2$ which in turn implies (i): see the proof of ``(iv) $\Rightarrow$ (i)". Moreover, by (4.5), both relations (4.26)-(4.27) imply the first limit (4.28).\par
(vi) $\Rightarrow$ (v). We already proved that (vi) implies (4.29)$_1$, i.e. (4.30)$_1$; moreover from (4.31)$_2$ we get
$$
(f(x)/\phi_2(x))'=a_1(\phi_1(x)/\phi_2(x))'-\left({\Phi(x)\over\phi_2(x)}\right)'\cdot \int\limits_x^{x_0} W\Phi^{-2} \overline F+
{W(x)\overline F(x)\over \phi_2(x)\Phi(x)}\,,\leqno(7.24)
$$
and we must estimate the last two terms on the right-hand side. We have
$$
\left({\Phi(x)\over\phi_2(x)}\right)'\cdot\int\limits_x^{x_0} W\Phi^{-2} \overline F\stackrel{\ by\ (7.9)\ }{=}\left({\Phi\over\phi_2}\right)'\cdot o\left({\phi_2\over\Phi}\right)=o\left({W(\phi_2,\Phi)\over\phi_2\Phi}\right)=\leqno(7.25)
$$
$$
\stackrel{\ by\ (7.1)\ }{=}o\left({W\over\phi_2\Phi}\right)\stackrel{\ by\ (7.7)\ }{=}o\left({W\over\phi_1\phi_2}\right)=o\left({\phi_2\over\phi_1}
\left({\phi_1\over\phi_2}\right)'\right)\,;
$$
and the last two passages show that the function $W\overline F(\phi_2\Phi)^{-1}$ satisfies the same asymptotic estimate as well. Substituting into (7.24) we get (4.30)$_2$.
\par
(v) $\Rightarrow$ (i). We may resort to representation (7.17) from which we get, using (4.7),
$$
(f(x)/\phi_2(x))'=a_1(\phi_1(x)/\phi_2(x))'-\phi_2(T)\left({\phi_1(x)\over\phi_2(x)}
\right)'\cdot \displaystyle{\int\limits_x^{x_0}} W\Phi^{-2} F^*
      +{W(x) F^*(x)\over\phi_2(x)\Phi(x)}.\leqno(7.26)$$
\par
For brevity we put
$$
\Omega(x):={\phi_2(x)\over\phi_1(x)}\left({\phi_1(x)\over\phi_2(x)}\right)'\equiv{-W(x)\over\phi_1(x)\phi_2(x)}\neq 0\quad \forall x\in I\,.\leqno(7.27)
$$
\par
Both (7.26) and (4.30)$_2$ imply
$$
{W(x)F^*(x)\over\phi_2(x)\Phi(x)}-\phi_2(T)\left({\phi_1(x)\over\phi_2(x)}\right)'\cdot \int\limits_x^{x_0}W\Phi^{-2} F^*=o(\Omega(x))\,.\leqno(7.28)
$$
\par
Now by (7.7):
$$
W(x)/\phi_2(x)\Phi(x)\sim {1\over \phi_2(T)}\,{W(x)\over \phi_1(x)\phi_2(x)}\equiv -{1\over \phi_2(T)}\,\Omega(x)\,;\leqno(7.29)
$$
and, by the already-proved theorem 4.4, condition (4.23) and (7.7) imply
$$
\int\limits_x^{x_0}W\Phi^{-2} F^*={\phi_2(x)\over\phi_1(x)}\, \left[-{l\over (\phi_2(T))^2} + o(1)\right]\,.\leqno(7.30)
$$
\par
Using (7.29) and (7.30) into (7.28) we easily get (4.26). Obviously (4.35) follows from (4.5) and representations (4.24), (4.25) hold true with $m=a_2$ and $l=\gamma$. The proof is complete. The claims in part II are contained in part I.
\hfill{$\Box$}
\par\vspace{15pt}
{\it Proof of theorem} 5.3.\ Representations (5.20) and (5.21) are nothing but (4.24) and (4.25); (5.22) and (5.24) are respectively obtained from (5.20) and (5.21) using (5.16) and (5.18); (5.23) follows from (4.20) using (5.13). In proving the estimates (5.25)-(5.26) we use in an essential way the constancy of sign of $(\phi_2/\phi_1)'$. From (5.20) we get
$$
\vert R(x)\vert\leq\vert\phi_1(x)\vert\cdot \left({\mathrel{\mathop{ess.sup.}_{x<t<x_0}}}\vert f^*_2(t)-a_2\vert\right)\cdot\left\vert\displaystyle{\int\limits_x^{x_0}}(\phi_2/\phi_1)'\right\vert=
\vert\phi_2(x)\vert\cdot\left(
{\mathrel{\mathop{ess.sup.}_{x<t<x_0}}}\vert f^*_2(t)-a_2\vert\right);$$
whereas from (5.21), using (7.9),we get
$$
\begin{array}{ll}
\vert R(x)\vert\kern-4pt&\leq\vert\Phi(x)\vert \left(\displaystyle{\mathrel{\mathop{ess.sup.}_{x<t<x_0}}}\vert F^*(t)-\gamma\vert\right)\cdot\left\vert\displaystyle{\int\limits_x^{x_0}} W(t) (\Phi(t))^{-2}dt\right\vert=\\ \\
&=\left\vert\displaystyle{\phi_2(x)\over\phi_2(T)}\right\vert\left(
\displaystyle{\mathrel{\mathop{ess.sup.}_{x<t<x_0}}}\vert F^*(t)-\gamma\vert\right)\,.
\end{array}
$$
\par
In a similar manner (5.27) and (5.28) are proved. Representation (5.29) follows from (5.20) by application of the mean-value theorem for improper integrals of type $\int_a^b\phi(t) f(t) dt$. If
$$\begin{cases}
\phi\in L^1 (a, b)\; \,; \; \phi\geq 0\,;\quad f\in C^0 (]a, b[)\,;\quad a,b\in \overline{\mathbb R}\,;\cr 
the\; limits\; f(a^+)\; and \; f(b^-)\; exist\; as\; finite\; numbers\,;\cr
\end{cases}\leqno(7.31)$$
then
$$
\left(\int\limits_a^b \phi\right)\cdot\left(\displaystyle{\mathrel{\mathop{inf}_{(a,b)}}}\; f\right)\leq \int\limits_a^b \phi \,f\leq\left(\int\limits_a^b \phi\right)\cdot\left(\displaystyle{\mathrel{\mathop{sup}_{(a,b)}}}\; f\right) \leqno(7.32)
$$
from whence
$$
\int\limits_a^b \phi f =f(\xi)\cdot \int\limits_a^b \phi\quad for\; a \; suitable\; \xi\in [a, b]\,.\leqno(7.33)
$$
\par
This last equality holds true for $\phi\leq 0$ as well. In the same way (5.30) is obtained from (5.21) using (7.9).
\hfill{$\Box$}
\par\vspace{15pt}
For the proof of theorems 6.1, 6.2 we use some known nontrivial results about generalized convex functions. Next lemma characterizes generalized convexity on an open interval when the underlying system comprises functions of class $C^1$ taking account of at least one of the endpoints.
\par\vspace{5pt}
\textbf{Lemma 7.4}\quad (Karlin and Studden  [10; chp XI; th. 2.1, p. 386]). \textsl{Let} $(\psi_1, \psi_2)$ \textsl{be an} $ECT$-\textsl{system on the interval} $[a, b[$ \textsl{of the explicit form}
$$
\psi_1(x):=w_1(x)\quad\,;\quad \psi_2(x):=w_1(x)\cdot\int\limits_a^x w_2(t) dt\,,\leqno(7.34)
$$
\textsl{wherein}
$$
w_1(x), w_2(x)>0\quad \forall x\in [a, b[\,;\quad 
w_1\in C^1 [a, b[\,;\quad w_2\in C^0 [a, b[\,.\leqno(7.35)
$$
\par
\textsl{Then} $f\in\mathcal C(\psi_1, \psi_2; ]a,b[)$ \textsl{iff all the following properties hold true}:
\par
(i)\quad $f$ \textsl{is continuous on } $]a, b[$;
\par
(ii)\quad $f$ \textsl{has right derivative $f'_R$ which is right-continuous on $]a, b[$ and a left derivative $f'_L$ which is left-continuous on} $]a, b[$;
\par
(iii)\quad \textsl{The function}
$$
\rho (x):={1\over w_2(x)}\cdot D_R (f(x)/w_1(x))\leqno(7.36)
$$
\textsl{is right-continuous and increasing on} $]a, b[$, $D_R$ \textsl{denoting the right derivative.}
\par\vspace{15pt}
\textsl{Remark}.
In the quoted reference the regularity assumptions (7.35) are assumed valid on a compact interval $[a, b]$ but this is immaterial for our thesis as $f\in \mathcal C(\psi_1,\psi_2; ]a, b[)$ iff $f\in\mathcal C(\psi_1,\psi_2; [\alpha, \beta])$ for each $\alpha, \beta$: $a<\alpha<\beta<b$.
\par
The next elementary lemma claims the invariance of the property of generalized convexity with respect to some systems to be used in the sequel.
\par\vspace{5pt}
\textbf{Lemma 7.5}\quad
\textsl{Let} $(\phi_1, \phi_2)$ \textsl{satisfy conditions} (6.6) \textsl{and define the following two} $ECT$-\textsl{systems on} $I$:
$$
\overline\psi_1(x):=\phi_1(x)\,;\; \overline\psi_2(x):=\phi_1(x)\cdot\int\limits_T^x (\phi_2(t)/\phi_1(t))' dt\equiv \phi_2(x)+\overline c\phi_1(x)\,;\leqno(7.37)
$$
$$
\begin{array}{lll}
\kern-1.9cm(7.38)&\overline{\overline \psi}_1(x):\kern-4pt&=\vert\phi_2(x)\vert\,;\; \overline{\overline{\psi}}_2(x):=-\vert\phi_2(x)\vert
\cdot \displaystyle{\int\limits_T^x} (\phi_1(t)/\phi_2(t))' dt\equiv \\ \\
&&\equiv(- sign \;\phi_2)\phi_1(x)+\overline{\overline c}\phi_2(x)\,.
\end{array}
$$\par
\textsl{Then} $f\in \mathcal C(\phi_1, \phi_2; I) \Leftrightarrow f\in \mathcal C(\overline\psi_1, \overline\psi_2; I)\Leftrightarrow f\in \mathcal C(\overline{\overline\psi}_1, \overline{\overline\psi}_2; I)$.
\par\vspace{15pt}
{\it Proof.}\ Trivial evaluations of the determinant (6.5) for each of the three systems $(\phi_1, \phi_2)$, $(\overline\psi_1, \overline\psi_2)$, $(\overline{\overline\psi}_1, \overline{\overline\psi_2})$ give the same value.
\par\vspace{5pt}
\textbf{Lemma 7.6}\quad (Characterizations of generalized convexity via the monotonicity of $f^*_1, f^*_2, F^*$). \textsl{If the pair} $(\phi_1, \phi_2)$ \textsl{satisfies conditions} (6.6) \textsl{and if} $\Phi$ \textsl{is defined by} (4.7) \textsl{then} $f\in\mathcal C (\phi_1, \phi_2; ]T, x_{{}_0}[)$ \textsl{iff all the following properties hold true}:
\par
(i)\quad $f$ \textsl{is continuous on } $]T, x_{{}_0}[$;
\par
(ii)\quad $f$ \textsl{has a right derivative which is right-continuous and a left derivative which is left-continuous on} $]T, x_{{}_0}[$;
\par
(iii)\quad \textsl{any of the three functions}
$$
\rho_1(x):={D_R (f(x)/\phi_1(x))\over (\phi_2(x)/\phi_1(x))'}\equiv {-W_R(f(x), \phi_1(x))
\over W(x)}\,,\leqno(7.39)
$$
$$
\rho_2(x):=(- sign\, \phi_2){D_R (f(x)/\phi_2(x))\over (\phi_1(x)/\phi_2(x))'}\equiv (- sign \,\phi_2)
 {W_R(f(x), \phi_2(x))\over W(x)}\,,\leqno(7.40)
$$
$$
\rho_3(x):=(sign \, \phi_2){W_R(\Phi(x), f(x))\over W(x)}\,,\leqno(7.41)
$$
\textsl{is right-continuous and increasing} ($\equiv$ \textsl{nondecreasing}) \textsl{on} $]T, x_{{}_0}[$. \textsl{Here} $W_R$ \textsl{denotes the Wronskian constructed with the right derivatives}.
\par
\textsl{Notice that at each point $t$ where $f$ is differentiable the three values} $(- sign\, \phi_2) \rho_2(t)$, $\rho_1(t)$ \textsl{and} $(sign\, \phi_2) \rho_3(t)$ \textsl{respectively coincide with} $f_1^*(t), f_2^*(t), F^*(t)$: the three geometric quantities upon which our theory rotates.
\par
\textsl{Under the stronger regularity conditions} (3.1)-(3.4) \textsl{and} $f\in AC^1 (I)$ \textsl{then} 
\par
$f\in\mathcal C(\phi_1, \phi_2; ]T, x_{{}_0}[)$ \textsl{iff} $Lf(x)\geq 0$ \textsl{a.e. on} $[T, x_{{}_0}[$.
\par\vspace{5pt}
{\it Proof.}\ The statements involving either $\rho_1$ or $\rho_2$ directly follow from lemmas 7.4, 7.5. Let us prove the statement involving $\rho_3$. Define
$$
\widetilde\psi_1 :=\vert\Phi(x)\vert\,\,;\,\widetilde\psi_2 :=\vert\Phi(x)\vert\cdot \int\limits_{T_0}^x\vert(\phi_2(t)/\Phi(t))'\vert dt\,,\quad x\in ]T, x_0[\,,\leqno(7.42)
$$
where $T_0$ is fixed in $]T, x_{{}_0}[$. From (4.8), (6.6)$_3$ and (7.7) we infer $sign\, \Phi\, = sign\, \phi_2$, and from (6.6)$_4$ and (7.11) we infer $sign\, (\phi_2/\Phi)'= sign\, W(\Phi, \phi_2)= sign\, \phi_2$. Hence we rewrite (7.42) as
$$
\widetilde\psi_1 :=(sign\, \phi_2) \Phi(x)\,\,;\,\widetilde\psi_2 :=\Phi(x) \cdot \int\limits_{T_0}^x\vert\phi_2(t)/\Phi(t))'dt\,,\quad x\in ]T_0, x[\,.\leqno(7.43)
$$
\par
By (4.7) we explicitly get
$$
\widetilde\psi_1(x)=(sign\, \phi_2) [\phi_2(T)\phi_1(x)-\phi_1(T)\phi_2(x)]\equiv\alpha\phi_1(x)+\beta\phi_2(x);
\leqno(7.44)_1
$$
$$
\widetilde\psi_2(x)=\left[\phi_2(x)-{\phi_2(T_0)\over\Phi(T_0)}\Phi(x)\right]=
\leqno(7.44)_2
$$
$$
=\left[ {-\phi_2(T_0)\phi_2(T)\over\Phi(T_0)} \phi_1(x)+\left(1+{\phi_2(T_0)\phi_1(T)\over\Phi(T_0)}\right)\phi_2(x)\right]\equiv
\overline\alpha\phi_1(x)+\overline\beta\phi_2(x);
$$
from whence
$$  
U\left(\begin{array}{ccc}
\tilde\psi_1, & \tilde\psi_2,& f \\
t_1, & t_2, & t_3 \end{array}\right)
=U\left(\begin{array}{ccc}
\alpha\phi_1+\beta\phi_2, &\overline\alpha\phi_1+\overline\beta\phi_2,& f\\
t_1, & t_2, & t_3\end{array}\right)=\leqno(7.45)
$$
$$
=U\left(\begin{array}{ccc}
\alpha\phi_1, &\overline\beta\phi_2,& f\\
t_1, & t_2, & t_3\end{array}\right)+U\left(\begin{array}{ccc}
\beta\phi_2, &\overline\alpha\phi_1,& f\\
t_1, & t_2, & t_3\end{array}\right)=(\alpha\overline\beta+\overline\alpha\beta)\cdot 
U\left(\begin{array}{ccc}
\phi_1, &\phi_2,& f\\
t_1, & t_2, & t_3\end{array}\right)\,,
$$
where, as trivially checked, $\alpha\overline\beta+\overline\alpha\beta=\vert\phi_2(T)\vert$.  This implies that
$$
f\in\mathcal C(\phi_1, \phi_2; ]T_0, x_0[ )\Leftrightarrow f\in\mathcal C(\tilde\psi_1, \tilde\psi_2; ]T_0, x_0[ )\,.
$$
\par
Applying lemma 7.4 to the system (7.43) we get our statement, referred to the interval $]T_0, x_{{}_0}[$, as the function (7.36) now becomes
$$
{(sign\, \phi_2)\cdot D_R (f(x)/\Phi(x))\over (sign\, \phi_2)\cdot(\phi_2(x)/\Phi(x))'}={W_R(\Phi(x), f(x))
\over W(\Phi(x), \phi_2(x))}\stackrel{(4.7)}{=}{W_R(\Phi(x), f(x))
\over \phi_2(T) W(x)}\,.\leqno(7.46)
$$
\par
The last statement concerning $Lf\geq 0$ directly follows from representation (5.5) and the increasing character of $f^*_2\equiv \rho_1$.
\hfill{$\Box$}
\par\vspace{5pt}
{\it Proof of theorem} 6.1.\ We report a proof of (i) as we do not have a reference for it in the literature. The monotonicity of the function (7.36) implies
$$\begin{cases}
\rho\quad \hbox{\textsl{continuous on }} ]T, x_0[\setminus\tilde N \  \hbox{\textsl{for a suitable countable set }} \tilde N\,,\cr
\rho\quad \hbox{\sl locally bounded on }\, ]T, x_0[\,.\cr\end{cases}\leqno(7.47)
$$
\par
By the continuity of $w_2$ we infer from (7.36) that the function $D_R(f(x)/w_1(x))$, which is defined everywhere on $]T, x_{{}_0}]$ as a finite number, enjoys of the same two properties listed in (7.47) for $\rho$. A known result on the regularity of derivatives implies that $(f(x)/w_1(x))'$ exists as a finite number on $]T, x_{{}_0}[\setminus\tilde N$ and is locally bounded on $]T, x_{{}_0}[$. A result in real analysis, Hewitt and Stromberg [9; exercise 18.41-(d), p. 299], now implies that $f/w_1\in AC (]T, x_{{}_0}[)$ from whence 
$f\equiv w_1 (f/w_1)\in AC (]T, x_{{}_0}[)$ \  as \  $w_1\in C^1(]T, x_{{}_0}[)\,.$
\par
(ii) is contained in lemma 7.6.
\par
(iii). Put $\tilde f:= f+a_1\phi_1+a_1\phi_2$; then $\tilde f\in\mathcal C (\phi_1, \phi_2; ]T, x_{{}_0}[ )$ and property (iii) of lemma 7.6 applied to $\tilde f$ implies that $(\tilde f/\phi_i)'$, $i=1,2$, is either $\equiv 0$ or strictly one-signed on the specified neighborhood and the statement is proved for $\tilde\phi=\phi_1, \phi_2$. If $\tilde\phi$ is any functions in $\mathcal F$, $\tilde\phi\not\equiv 0$, then we kow that it has at most one zero on $]T, x_{{}_0}[$ hence it is $\neq 0$ on suitable deleted neighborhoods of $T$ and of $x_{{}_0}$. For any such neighborhood it is obviously possible to choose another function $\tilde{\tilde\phi\,}$ in $\mathcal F$ such that either the pair $(\tilde\phi, \tilde{\tilde\phi\,})$ or $(\tilde{\tilde\phi\,}, \tilde\phi)$ may play the role that $(\phi_1, \phi_2)$ played in lemmas 7.5 and 7.6 except possibly for some sign. Moreover either $f$ or $-f$ belongs to the cone $\mathcal C (\tilde\phi, \tilde{\tilde\phi\,})$ in the chosen neighborhood and lemma 7.6, referred to this new context, implies that $(\tilde f/\tilde\phi)'$ is either $\equiv 0$ or strictly one-signed on suitable neighborhoods of the endpoints.
\par
(iv) follows from the monotonicity of the two quotients on the two sides of (6.7) and from L'Hospital's rule when one writes $f/\phi_1=(f/\phi_2)/(\phi_1/\phi_2)$. The functions $\phi_1, \phi_2$ cannot be interchanged because L'Hospital's rule only works when the denominator tends to $\pm\infty$ if no information is available on the numerator.
\hfill{$\Box$}
\par\vspace{5pt}
{\it Proof of theorem} 6.2.\ (i). If $f=O(\phi_1)$ then property (iii) in theorem 6.1 implies the existence of a finite $\lim_{x\to x_0} f/\phi_1$, and the second asymptotic relation follows from (6.7).
\par
(ii). The statement involving (6.9)-(6.10) are obvious by the monotonicity of $F^*$ and $f^*_2$. Now we show that (6.11) implies (4.29)$_{1,2}$. Relation (4.29)$_1$ follows from the ultimate monotonicity of $(f-a_1\phi_1)/\phi_2$, whereas (4.29)$_2$ follows from the following application of L'Hospital's rule
$$
a_2=\lim_{x\to x_0} {f(x)-a_1\phi_1(x)\over \phi_2(x)}\equiv \lim_{x\to x_0} {{f(x)\over \phi_1(x)}-a_1\over \phi_2(x)/\phi_1(x)}=
\lim_{x\to x_0} {(f(x)/\phi_1(x))'\over (\phi_2(x)/\phi_1(x))'}\,,
\leqno(7.48)
$$
as the last quotient is ultimately monotonic. 
\par
(iii). Inequalities (6.12), (6.13), (6.14) follow from the types of monotonicity of the involved functions,
whereas (6.15) follows from, say, representation (4.24) with $m=a_2$, and from (6.14). The last assertion about $R(x)$ again follows from (4.24), or equivalently from (4.25), and the following two facts: the function $R(x)/\phi_1(x)$ tends to zero, as $x\to x_{{}_0}$, and is either constant or strictly monotonic on some left neighborhood of $x_{{}_0}$.
\hfill{$\Box$}
\par\vspace{20pt}
\centerline{\bf 8. Example: the special case of powers.}
\par\vspace{10pt}
The case discussed in this section also serves as an illustration for the theory developed in [6; 7].
To fix the ideas we suppose $f\in AC^1 ]0, +\infty)$ and study the validity of asymptotic expansions of type
$$
\begin{array}{ll}
\kern-2cm(8.1)\qquad &f(x)=a_1 x^{\alpha_1}+a_2 x^{\alpha_2}+o(x^{\alpha_2})\,,\quad x\to +\infty \\ \\
\kern-2cm(8.2)\qquad &f(x)=a_2 x^{\alpha_2}+a_1 x^{\alpha_1}+o(x^{\alpha_1})\,,\quad x\to 0^+ \\ 
\end{array} \qquad (\alpha_1>\alpha_2),
$$
where $\alpha_1, \alpha_2$ are arbitrary real numbers. The associated (Euler) differential operator is
$$
L_{\alpha_1, \alpha_2} [u] := u ''+(1-\alpha_1-\alpha_2) x^{-1} u'+\alpha_1\alpha_2 x^{-2} u ''\,,\quad x>0\,.\leqno(8.3)
$$
\par
As convenient canonical factorizations we may use the following
$$
L_{\alpha_1, \alpha_2} [u] \equiv x^{\alpha_1-1} 
\left[ x^{\alpha_2-\alpha_1+1}\left(x^{-\alpha_2}u\right)'\right]'
{\nearrow\atop{\searrow }}
\begin{array}{ll}
\hbox{type (I) at} +\infty \\ \\
\hbox{type (II) at}\;\,0^+
\end{array}\quad\,; \leqno(8.4)
$$
$$
L_{\alpha_1, \alpha_2} [u] \equiv x^{\alpha_2-1} 
\left[ x^{\alpha_1-\alpha_2+1}\left(x^{-\alpha_1}u\right)'\right]'
{\nearrow\atop{\searrow }}
\begin{array}{ll}
\hbox{type (II) at} +\infty \\ \\
\hbox{type (I) at}\;\, 0^+
\end{array}\quad\,. \leqno(8.5)
$$
\par
The simple choice $T=T_0=1$ in the formulas of \S 4 is an admissible and convenient one for all values of $\alpha_i$
\par
\underline{Case} : $x\to +\infty$. The basic quantities of our theory are:
$$\begin{cases}
\phi_1(x) := x^{\alpha_1}\,\,;\; \phi_2(x) :=x^{\alpha_2}\,;\; W(x)=(\alpha_2-\alpha_1) x^{\alpha_1+\alpha_2-1}\quad (x>0)\,;\cr \cr
\Phi(x)=x^{\alpha_1}-x^{\alpha_2}\,.\cr\end{cases}\leqno(8.6)
$$
$$
f^*_1(x)={-x^{\alpha_2}f '(x)+\alpha_2 x^{\alpha_2-1}f(x)\over (\alpha_2-\alpha_1) x^{\alpha_1+\alpha_2-1}}=
{-x f '(x)+\alpha_2 f(x)\over (\alpha_2-\alpha_1) x^{\alpha_1}}=\leqno(8.7)
$$
$$
\stackrel{\ by\ (5.4)\ }{=} \overline c_1+{1\over \alpha_1-\alpha_2}\int\limits_1^x t^{1-\alpha_1} \cdot L_{\alpha_1, \alpha_2} [f(t)] dt\ ,\quad \overline c_1=\alpha_2 f(1)-f '(1)\,;
$$
$$
f^*_2(x)={x^{\alpha_1}f '(x)-\alpha_1 x^{\alpha_1-1}f(x)\over (\alpha_2-\alpha_1) x^{\alpha_1+\alpha_2-1}}=
{x f '(x)+\alpha_1 f(x)\over (\alpha_2-\alpha_1) x^{\alpha_2}}=\leqno(8.8)
$$
$$
\stackrel{\ by\ (5.5)\ }{=} \overline c_2+{1\over \alpha_2-\alpha_1}\int\limits_1^x t^{1-\alpha_2} \cdot L_{\alpha_1, \alpha_2} [f(t)] dt\ ,\quad \overline c_2=f '(1)-\alpha_1 f(1)\,;
$$
$$
F^*(x)={W(\Phi(x), f(x)) \over W(x)}=
{(x^{\alpha_1}-x^{\alpha_2}) f '(x)- (\alpha_1 x^{\alpha_1-1}-\alpha_2 x^{\alpha_2-1}) f(x)
\over (\alpha_2-\alpha_1) x^{\alpha_1+\alpha_2-1}}=\leqno(8.9)
$$
$$
\stackrel{\ by\ (5.6)\ }{=}f(1)+{1\over \alpha_2-\alpha_1}\int\limits_1^x {t^{\alpha_1}-t^{\alpha_2}\over t^{\alpha_1+\alpha_2-1}} 
\cdot L_{\alpha_1, \alpha_2} [f(t)] dt\,.
$$
\par
Specializing our theory we get the following results.
\par\vspace{5pt}\textbf{Proposition 8.1.}\quad Part I. \textsl{The following are equivalent properties:}
\par
(i) \textsl{The pair of asymptotic expansions}
$$\begin{cases}
f(x)= a_1 x^{\alpha_1}+a_2 x^{\alpha_2}+o(x^{\alpha_2})\cr\cr
f '(x)= a_1 \alpha_1 x^{\alpha_1-1}+a_2\alpha_2 x^{\alpha_2-1}+o(x^{\alpha_2-1})\cr\end{cases}\quad\,, x\to +\infty\,;\leqno(8.10)
$$
\par
(ii) \textsl{The pair of asymptotic expansions}
$$\begin{cases}
f(x)= a_1 x^{\alpha_1}+a_2 x^{\alpha_2}+o(x^{\alpha_2})\cr\cr
(x^{-\alpha_1}f(x))'= a_2 (\alpha_2-\alpha_1) x^{\alpha_2-\alpha_1-1}+o(x^{\alpha_2-\alpha_1-1})\cr\end{cases}\quad\,, x\to +\infty\,;\leqno(8.11)
$$
\par
(iii) \textsl{The pair of asymptotic expansions}
$$\begin{cases}
f(x)= a_1 x^{\alpha_1}+a_2 x^{\alpha_2}+o(x^{\alpha_2})\cr\cr
(x^{-\alpha_2}f(x)) '= a_1 (\alpha_1-\alpha_2) x^{\alpha_1-\alpha_2-1}+o(x^{-1})\cr\end{cases}\quad\,, x\to +\infty\,;\leqno(8.12)
$$
\par
(iv) \textsl{The improper integral}
$$
\int\limits_1^{+\infty} t^{1-\alpha_2}\cdot L_{\alpha_1, \alpha_2} [f(t)] dt\equiv \int\limits_1^{+\infty} [t^{\alpha_1-\alpha_2+1}(t^{-\alpha_1} f(t))']' dt\;\; converges\,. \leqno(8.13)
$$
\par
\textsl{And to this list we may add the geometric properties in theorem 4.5 concerning the limits of} $f_1^*, f_2^*, F^*$.
\par
Part II. \textsl{Whenever properties in part I hold true then}
$$
a_1=\alpha_2 f(1)-f '(1)+{1\over \alpha_1-\alpha_2}
\int\limits_1^{+\infty} [t^{\alpha_2-\alpha_1+1}(t^{-\alpha_2} f(t))']' dt\,; \leqno(8.14)
$$
$$
a_2=-\alpha_1 f(1)+f '(1)+{1\over \alpha_2-\alpha_1}
\int\limits_1^{+\infty} [t^{\alpha_1-\alpha_2+1}(t^{-\alpha_1} f(t))']' dt\,; \leqno(8.15)
$$
\textsl{and we have the following representations}
$$
f(x)=a_1 x^{\alpha_1}+ a_2 x^{\alpha_2}+x^{\alpha_1}\cdot
\int\limits_x^{+\infty} t^{\alpha_2-\alpha_1-1} dt \int\limits_t^{+\infty} s^{1-\alpha_2}\cdot 
L_{\alpha_1,\alpha_2} [f(s)] ds\,; \leqno(8.16)
$$
$$
f(x)=a_1 x^{\alpha_1}+ a_2 x^{\alpha_2}+x^{\alpha_2}\cdot
\int\limits_x^{+\infty} t^{\alpha_1-\alpha_2-1} dt \int\limits_t^{+\infty} s^{1-\alpha_1}\cdot 
L_{\alpha_1,\alpha_2} [f(s)] ds\,. \leqno(8.17)
$$
\par
Part III. \textsl{In the special case wherein the quantity} $L_{\alpha_1, \alpha_2} [f(x)]$ \textsl{is one-signed} ($\geq 0$ \textsl{or} $\leq 0$) \textsl{for all $x$ large enough then to all the equivalent properties listed in part I the following can be added:}
$$
f(x)=a_1 x^{\alpha_1}+ O(x^{\alpha_2})\quad\,,\quad x\to +\infty \,;\leqno(8.18)
$$
$$
f(x)=a_1 x^{\alpha_1}+a_2 x^{\alpha_2}+o(x^{\alpha_2})\,,\quad x\to +\infty\,.\leqno(8.19)
$$
\par
\textsl{Here the import is that the sole relation} (8.18) \textsl{automatically implies the pair} (8.10).
\par
Part IV. \textsl{The following are equivalent properties}
\par
(v) \textsl{The pair of asymptotic expansions}
$$\begin{cases}
f(x)=a_1 x^{\alpha_1}+a_2 x^{\alpha_2}+o(x^{\alpha_2})\cr\cr
f '(x)=a_1 \alpha_1 x^{\alpha_1-1}+o(x^{\alpha_1-1})\cr\end{cases}\quad\,,\quad x\to +\infty\,;\leqno(8.20)$$
\par
(vi) \textsl{The pair of asymptotic expansions}
$$\begin{cases}
f(x)=a_1 x^{\alpha_1}+a_2 x^{\alpha_2}+o(x^{\alpha_2})\cr\cr
(x^{-\alpha_2} f(x)) '=a_1 (\alpha_1-\alpha_2)x^{\alpha_1-\alpha_2-1}
+o(x^{\alpha_1-\alpha_2-1})\cr\end{cases}\quad\,,\quad x\to +\infty\,;\leqno(8.21)
$$
\par
(vii) \textsl{The improper integral}
$$
\int\limits_1^{+\infty} t^{\alpha_1-\alpha_2-1} dt \int\limits_t^{+\infty} s^{1-\alpha_1}\cdot 
L_{\alpha_1,\alpha_2} [f(s)] ds\; converges\,. \leqno(8.22)
$$
\par\vspace{1pt}
The two equivalences ``(8.10) $\Leftrightarrow$ (8.11)" and ``(8.20) $\Leftrightarrow$ (8.21)" are not contained in the theory developed in this paper
but are simple algebraic facts that can be directly checked for any numbers $\alpha_1, \alpha_2\in\mathbb R$. The corresponding proofs for $n$-term expansions in real powers are to be found in [7; lemmas 7.3 and 7.4]. In this case the factorizational approach gives characterizations of standard differentiation of an asymptotic expansion, i.e. differentiation obtained by the application of the operator $d/dx$.
\par
This fact can be extended to a larger class of asymptotic expansions using the concept of regular variation but it will not be investigated here.
\par
For the elementary case of asymptotic straight lines, i.e. $\alpha_1=1$ and $\alpha_2=0$, we have the characterizations
$$\begin{cases}
f(x)= a_1 x+a_2+o(1)\cr\cr
f '(x)= a_1+o(1)\cr\end{cases}\quad\,, x\to +\infty\quad\Leftrightarrow\;
\int\limits_1^{+\infty} dt \int\limits_t^{+\infty} f ''(s) ds\; converges\,;\leqno(8.23)
$$
$$\begin{cases}
f(x)= a_1 x+a_2+o(1)\cr\cr
f '(x)= a_1+o(x^{-1})\cr\end{cases}\quad\,, x\to +\infty\quad\Leftrightarrow\;
\int\limits_1^{+\infty} t f ''(t) dt\; converges\,.\leqno(8.24)
$$
\par
The pair of asymptotic relations in (8.23) may be labelled by the locution ``\textsl{the graph of $f$ admits of the straight line $y=a_1 x+a_2$ as a first-order asymptote at} $+\infty$". The pair in (8.24) states the fact that the straight line $y=a_1 x+a_2$ is ``\textsl{the limit tangent} ($\equiv$ \textsl{asymptotic tangent}) \textsl{at} $+\infty$".
\par\vspace{5pt}
\underline{Case}: $x\to 0^+$. The basic quantities are the same as in the foregoing case with the roles of $\alpha_1, \alpha_2$ interchanged; and each integral of type, say $\int_1^{+\infty}\dots$ must be replaced by the integral $\int_{\to 0}^{1}\dots$ of the same quantity, wherein the endpoint ``$0$" is the sole possible singularity. We leave to the reader the complete formulation of the corresponding version of proposition 8.1 mentioning only the two main equivalences.
\par
\textsl{The pair of expansions}
$$\begin{cases}
f(x)=a_2 x^{\alpha_2}+a_1 x^{\alpha_1}+o(x^{\alpha_1})\cr\cr
f '(x)=a_2 \alpha_2 x^{\alpha_2-1}+a_1 \alpha_1 x^{\alpha_1-1}
+o(x^{\alpha_1-1})\cr\end{cases}\quad\,,\quad x\to 0^+\,,(\alpha_2<\alpha_1)\,,\leqno(8.25)$$
\textsl{holds true iff}
$$
\int\limits_{\to 0}^1 t^{1-\alpha_1}\cdot L_{\alpha_1,\alpha_2} [f(t)] dt\equiv \int\limits_{\to 0}^1 [t^{\alpha_2-\alpha_1+1}(t^{-\alpha_2} f(t))']' dt\;\, converges\,; \leqno(8.26)
$$
\textsl{whereas the pair of expansions}
$$\begin{cases}
f(x)=a_2 x^{\alpha_2}+a_1 x^{\alpha_1}+o(x^{\alpha_1})\cr\cr
f '(x)=a_2 \alpha_2 x^{\alpha_2-1}+o(x^{\alpha_2-1})\cr\end{cases}\quad\,,\quad x\to 0^+\,,(\alpha_2<\alpha_1)\,,\leqno(8.27)$$
\textsl{holds true iff}
$$
\int\limits_{\to 0}^1 t^{\alpha_2-\alpha_1-1} dt \int\limits_0^t s^{1-\alpha_2}
\cdot L_{\alpha_1,\alpha_2} [f(s)] ds\;\, converges\,. \leqno(8.28)
$$
\par
For $\alpha_2=0$ and $\alpha_1=1$ (8.25) reduces to
$$\begin{cases}
f(x)=a_2+a_1 x+o(x)\cr\cr
f '(x)=a_1+o(1)\cr\end{cases}\quad\,,\quad x\to 0^+\,,\leqno(8.29)
$$
which is obviously equivalent to
$$
f '(x)=a_1+o(1)\,,\quad x\to 0^+\,,\leqno(8.30)
$$
that is to say, to the existence of a finite limit: $\lim_{x\to 0^+} f '(x)$.
Condition (8.26) reduces to the convergence of $\int_{\to 0}^1 f ''(t) dt$, and this condition is equivalent to (8.30) under our present assumption $f\in AC^1 ]0, 1]$. 
\par
This is just the simple  technical idea underlying our factorizational theory together with the theory of canonical factorizations which yield the means for applying the simple idea to general expansions.
\par
In closing this paper we present a figure illustrating the concept of limit tangent curve, as characterized in theorem 4.5, for a generalized convex function.
\begin{figure}[htbp]
\begin{center}
\includegraphics[width=14cm]{grafico.jpg}
\end{center}
	\label{}
\end{figure}
\par\noindent
The figure refers to the following contingency:
\begin{enumerate}
	\item 
	$\phi_1, \phi_2, f\in C^2 [T, +\infty)$; $\phi_2>0$; $f\in\mathcal C(\phi_1, \phi_2; [T, +\infty))$, hence $F^*$ is increasing by theorem 6.1-(ii).
	\item
	Each dotted curve, save the uppermost, has a first-order contact with the graph of $f$ at a point $x_i$ and has equation
	$$
	y=f^*_1(x_i) \phi_1(x)+f^*_2(x_i) \phi_2(x)\,.
	$$
	\item 
	$F^*(x_i)$ is the contact indicatrix of order one at the point $x_i$ with respect to the family $\mathcal F := span (\phi_1, \phi_2)$ and to the line $x=T$.
	\item
	The $\lim_{x\to +\infty} F^*(x)=\gamma$ exists in $\mathbb R$.
\end{enumerate}

From these facts the following follow:
\par
(i) The two limits
$$
\lim_{x\to +\infty} f_1^*(x)\equiv a_1\,,\;\lim_{x\to +\infty} f_2^*(x)\equiv a_2\,,
$$
exist in $\mathbb R$ and are linked to $\gamma$ by relation 
$\gamma=a_1\phi_1(T)+a_2\phi_2(T)\,.$
\par
(ii) The uppermost dotted curve, whose equation is
$y=a_1\phi_1(x)+a_2(\phi_2(x)\,,$
is by definition the limit tangent curve to the graph of $f$ with respect to the family $\mathcal F$ as $x\to +\infty$.
\par
(iii) The asymptotic relations hold true:
$$\begin{cases}
f(x)=a_1\phi_1(x)+a_2\phi_2(x)+o(\phi_2(x))\,, x\to +\infty\,,\\ \\\displaystyle
\left(\displaystyle{f(x)\over\phi_1(x)}\right)'=a_2\left({\phi_2(x)\over\phi_1(x)}\right)'+o\left[\left(\displaystyle{\phi_2(x)\over\phi_1(x)}\right)'\right]\,,\quad x\to +\infty\,.\end{cases}$$

\vspace{15pt}
\begin{center}\bf 9.\  Formal differentiation of a two-term asymptotic expansion:\\ a Tauberian result of  interpolatory type\end{center}

\vspace{10pt}
Here we examine the classical problem of lookig for conditions under which an expansion
$$
f(x)=a_1\phi_1(x)+a_2\phi_2(x)+o(\phi_2(x))\,, x\to x_0,\leqno(9.1)$$
implies
$$
f'(x)=a_1\phi_1'(x)+a_2\phi_2'(x)+o(\phi_2'(x))\,, x\to x_0.\leqno(9.2)$$

The reder is referred to the introduction in [6] highlighting the inherent differences between the problems of differentiating in some formal sense an asymptotic relation 
\par\noindent
 $f(x)=\phi(x)+o\big(\phi(x)\big)$ or an asymptotic expansion with two meaningful terms.

For our problem the first remark is that theorems in \S\S4,5 show that the differentiated relation (9.2) is in general no good match for (9.1). As a simple example take the function
$$
g_1(x):=e^x+x+\sin x=e^x+x+o(x),\ x\to+\infty,\leqno(9.3)$$
defined, say, on $[2,+\infty)$ for which we have:
$$
g_1^{(k)}(x)=e^x+O(1),\ x\to+\infty;\ k\ge1.\leqno(9.4)$$

Hence the asymptotic expansion (9.3) is not formally differentiable if standard derivatives are used; however we  have:
$$
g_1(x)/e^x=1+xe^{-x}+\sin x\cdot e^{-x}=1+xe^{-x}[1+o(1)],\ x\to+\infty,\leqno(9.5)$$
$$
(g_1(x)/e^x)'=(1-x)e^{-x}+(\cos x-\sin x)e^{-x}=(1-x)e^{-x}[1+o(1)],\ x\to+\infty. \leqno(9.6)$$

Putting $\phi_1(x):=e^x,\ \phi_2(x):=x$ relation (9.6) is just (4.29)$_2$. The operator $L$ associated to the pair $(e^x,x)$ is 
$$
L[u]:=u''+{x\over1-x}u'-{1\over1-x}u\ \ on\ ]1,+\infty)\ ,\leqno(9.7)$$
and the integral in (5.15) becomes
$$
\int^{+\infty}{e^x\over (1-x)e^x}L[g(x)]dx=\int^{+\infty}\left[{-2\sin x+
x(\sin x+\cos x)\over(1-x)^2}\right]\ dx,\leqno(9.8)$$
which is convergent by Abel's test, according to Theorem 5.2-(II). A second simple example, left to the reader, is that of 
$$
g_2(x):=x+e^{-x}+x^{-1}\sin x\cdot e^{-x}=x+e^{-x}[1+o(1)],\ x\to+\infty,\leqno(9.9)$$
which is not formally differentiable; here $g_2$ has the properties in Theorem 4.3 but not those in Theorem 4.5.

Hence our theory puts in evidence the fact that the spontaneous choice of the operator $d/dx$ not always is the right choice when formally differentiating an asymptotic expansion with at least two meaningful terms. In the framework of our theory it is one of the two operators $L_1[u]:=(u/\phi_1)',\ L_2[u]:=(u/\phi_2)'$ which works well and the results in \S5 characterize the pair (4.19) and the pair (4.29)$_{1,2}$. These results, being characterizations, completely settle the problem but the classical formulation of the interpolatory approach was a bit different as will be clearly shown by the case of real powers. By Proposition 8.1 the pair of relations
$$
f(x)= a_1 x^{\alpha_1}+a_2 x^{\alpha_2}+o(x^{\alpha_2}),\ x\to+\infty,\ (\alpha_1>\alpha_2),\leqno(9.10)$$
$$
f '(x)= a_1 \alpha_1 x^{\alpha_1-1}+a_2\alpha_2 x^{\alpha_2-1}+o(x^{\alpha_2-1}),\ x\to+\infty,\ (\alpha_1>\alpha_2),\leqno(9.11)$$
is characterized by the integral condition in (8.13) with $L_{\alpha_1,\alpha_2}$ defined in (8.3). Hence knowing (9.10) to be true, relation (9.11) holds true iff (8.13) is satisfied. This is certainly the case if
$$
L_{\alpha_1,\alpha_2}=O\big(x^{\alpha_2-2-\epsilon}\big),\ x\to+\infty ,\ 
for\ some\ \epsilon>0;\leqno(9.12)$$
but in the study of the $n$-body problem, for instance, it is of interest to grant (9.11) under the weaker condition 
$$
L_{\alpha_1,\alpha_2}=O\big(x^{\alpha_2-2}\big),\ x\to+\infty,\leqno(9.13)$$
and this is no elementary question.

The problem can be posed in a technically different way by putting
$$
R(x):=f(x)- a_1 x^{\alpha_1}-a_2 x^{\alpha_2}$$
and inferring $R'(x)=o(x^{\alpha_2-1})$ from both relations 
$$
R(x)=o(x^{\alpha_2}),\ \ R''(x)=O(x^{\alpha_2-2}),\ x\to+\infty .$$

This inference is known to be true: an ``$\epsilon-\delta$''-proof may be found, e.g., in Boas [15] for a twice-differentiable function and in Saari [16] for an $f\in AC^1[T,+\infty).$

Now condition $R''(x)=O(x^{\alpha_2-2})$ is algebraically natural whereas condition (9.13) naturally follows from our factorizatinal theory. But in the case of a general asymptotic expansion there is no algebraic evidence and it is our theory that leads to formulate the appropriate

\vspace{5pt}
 {\bf Conjecture on the formal differentiation of a general two-term expansion from the classical interpolatory standpoint.}\ {\it Let a function $f\in AC^1[T,x_0[$ admit of an 
asymptotic expansion $(9.1)$ under our strenghtened basic asssumptions $(3.1)$ and $(3.4)$.
It follows from Theorems $4.5$ and {\rm 5.2-(II)} that $(9.1)$ is formally differentiable in the sense of $(4.29)_2$ iff
$$
\int\limits_T^{x_0} \phi_1(t) (W(t))^{-1}\cdot L[f(t)]dt\quad converges\ .\leqno(9.14)$$

We now suggest the following heuristic considerations. It is implcit in the use of an 
asymptotic expansion $(9.1)$ that we are measuring our quantities by means of the given functions $\phi_1, \phi_2$ and of their ratios for which we know that $\phi_2/\phi_1=o(1)$ and that $\int^{x_0}(\phi_2/\phi_1)'$\ converges. Hence in this context, if the integrand in $(9.14)$ is $O((\phi_2(x)/\phi_1(x))')$ i.e.
 $$
 L[f(x)]=O\left(\big[(\phi_2(x)/\phi_1(x))'\big]^2\cdot\phi_1(x)\right)\equiv
 O\left((W(x))^2(\phi_1(x))^{-3}\right),\ x\to x_0,\leqno(9.15)$$
  then $(9.14)$ is satisfied; but if this integrand satisfies
$$
\phi_1(x) (W(x))^{-1}\cdot L[f(x)]=\left.\begin{cases}O\\o\end{cases}\hspace{-11pt}\right\}
\left({\phi_1(x)\over \phi_2(x)}\left({\phi_2(x)\over \phi_1(x)}\right)'\right)\ ,\ \  x\to x_0,\leqno(9.16)$$
that is to say
$$
L[f(x)]=O\left(\big[(\phi_2(x)/\phi_1(x))'\big]^2 (\phi_1(x)^2(\phi_2(x))^{-1}\right)\equiv
 O\left((W(x))^2(\phi_1(x))^{-2}(\phi_2(x))^{-1}\right), x\to x_0,\leqno(9.17)$$
then $(9.14)$ is not automatically granted. This is precisely the Tauberian condition we wish to investigate conjecturing that it is sufficient for the inference 
``$(9.1)\implies (4.29)_2$}''.  In the case of powers $\phi_i\equiv x^{\alpha_i},\ \alpha_1>\alpha_2$, and as $x\to+\infty$, we have 
$$\begin{cases}
(\phi_2(x)/\phi_1(x))'=(\alpha_2-\alpha_1)x^{\alpha_2-\alpha_1-1},\\
(\phi_1(x)/\phi_2(x))\cdot(\phi_2(x)/\phi_1(x))'=(\alpha_2-\alpha_1)x^{-1}.\end{cases}\leqno(9.18)$$

\vspace{5pt}
We present here one of the possible
 results to show the usefulness of canonical factorizations in this context; the proof does not follow classical patterns but is based on an interplay between the two types of factorizations.

\vspace{10pt}
{\bf Theorem 9.1.}\ {\it Hypotheses: {\rm (i)} the strenghtened basic asssumptions $(3.1)$ and $(3.4)$;  {\rm (ii)} $f\in AC^1[T,x_0[$ and the expansion $(9.1)$ i.e. $(4.29)_1$;  {\rm (iii)} $L$ is the operator defined by {\rm (3.5)-(3.6)}. Thesis:

{\rm (I)}\ If $(9.17)$ is satisfied with ``$O$'' replaced by ``$o$'' then $(4.29)_2$ holds true.

{\rm (II)}\ If $(9.17)$ is satisfied then, in general, only the weaker relation holds true:
$$
\big(f(x)/\phi_1(x)\big)'=O\big((\phi_2(x)/\phi_1(x))'\big),\ x\to x_o.\leqno(9.19)$$

But if the ratio $\phi:=\phi_2/\phi_1$ satisfies the additional conditions that $\phi'$ is strictly one-signed on a neighborhood of $x_0$ and
$$\begin{cases}
(x-x_0)\phi'(x)\asymp\phi(x)\\ (x-x_0)\phi''(x)\asymp\phi'(x)
\end{cases},\ x\to x_0\ \ (if\ x_0\in {\mathbb R}),\leqno(9.20)$$
or
$$\begin{cases}
x\phi'(x)\asymp\phi(x)\\ x\phi''(x)\asymp\phi'(x)
\end{cases},\ x\to +\infty\ \ (if\ x_0=+\infty),\leqno(9.21)$$
then $(4.29)_2$ holds true.}

\vspace{5pt}
Notation $g_1(x)\asymp g_2(x),\ x\to x_0$, denotes the validity of both relations
$$
g_1(x)=O\big(g_2(x)\big);\ \ g_2(x)=O\big(g_1(x)\big),\ x\to x_0.\leqno(9.22)$$

Conditions (9.20) and (9.21) respectively imply
$$
\big(\phi'(x)\big)^2\asymp \phi(x)\phi''(x),\ \ \begin{cases}x\to x_0\\x\to+\infty
\end{cases},\leqno(9.23)$$
which will be essential in the proof.

\vspace{5pt}
{\it Proof.}\ (I)\ Put
$$
M_1(x):=\phi_1(x) (W(x))^{-1}\cdot L[f(x)]; \ \ M_2(x):=\phi_2(x) (W(x))^{-1}\cdot L[f(x)].\leqno(9.24)$$

Condition (9.17) with ``$O$'' replaced by ``$o$'' is equivalent to each one of the following:
$$
M_1(x):=o\left({\phi_1(x)\over \phi_2(x)}\left({\phi_2(x)\over \phi_1(x)}\right)'\right)
\equiv o\left({W(x)\over \phi_1(x)\phi_2(x)}\right) ,\  x\to x_0;\leqno(9.25)$$
$$
M_2(x):=o\left(\left({\phi_2(x)\over \phi_1(x)}\right)'\right)
\equiv o\left({W(x)\over \ (\phi_1(x))^2}\right) ,\  x\to x_0.\leqno(9.26)$$

Condition (9.25) by itself does grant neither (5.15) nor (4.14) but (9.26) grants (5.12) and so we may rewrite representation (5.9) in the form:
$$\begin{cases}\displaystyle
f(x)=c_1\phi_1(x)+c_2\phi_2(x)+\phi_2(x)\int\limits_T^x\big(\phi_1(t)/\phi_2(t)\big)'
dt\int\limits_t^{x_0}M_2(s)ds=\\ \\ \displaystyle
=c_1\phi_1(x)+c_2\phi_2(x)+\phi_2(x)\cdot o\left({\phi_1(x)\over \phi_2(x)}\right)=
c_1\phi_1(x)+o\big(\phi_1(x)\big),\  x\to x_0.
\end{cases}\leqno(9.27)$$

From (9.27) and (9.1) we get both $c_1=a_1$ and 
$$\begin{cases}\displaystyle
\int\limits_T^x\big(\phi_1(t)/\phi_2(t)\big)'dt\int\limits_t^{x_0}M_2(s)ds=
\big[f(x)-a_1\phi_1(x)-c_2\phi_2(x)\big]\big/\phi_2(x)=\\\displaystyle
\hspace{5pt}=\big[(a_2-c_2)\phi_2(x)+o\big(\phi_2(x)\big)\big]\big/\phi_2(x)=(a_2-c_2)+o(1),\  x\to x_0.\end{cases}\leqno(9.28)$$

Hence (5.14) holds true and, instead of (9.27) we may use representation
$$
f(x)=a_1\phi_1(x)+a_2\phi_2(x)-\phi_2(x)\int\limits_x^{x_0}\big(\phi_1(t)/\phi_2(t)
\big)'dt\int\limits_t^{x_0}M_2(s)ds,\ x\in [T,x_0[,\leqno(9.29)$$
where $a_1,a_2$ are the same coefficients as in (4.29)$_1$. Moreover Theorem 5.2-(I) implies the formally-differentiated expansion appearing in (4.19) but we want to prove 
the stronger relation (4.29)$_2$. From (9.29) we get
$$\begin{cases}\displaystyle
(f/\phi_1)'=a_2(\phi_2/\phi_1)'-(\phi_2/\phi_1)'\cdot \int\limits_x^{x_0}\big(\phi_1/\phi_2\big)'dt\int\limits_t^{x_0}M_2(s)ds+
\\ \displaystyle
+(\phi_2/\phi_1)\cdot (\phi_2/\phi_1)'\cdot\int\limits_x^{x_0}M_2(s)ds,
\end{cases}\leqno(9.30)$$
where, by (9.26), the last term in the right-hand side satisfies the estimate:
$$\begin{cases}\displaystyle
{\phi_2(x)\over \phi_1(x)}\left({\phi_1(x)\over \phi_2(x)}\right)'\cdot\int\limits_x^{x_0}M_2(s)ds\equiv -{\phi_1(x)\over \phi_2(x)}\left({\phi_2(x)\over \phi_1(x)}\right)'\cdot\int\limits_x^{x_0}M_2(s)ds=
\\ \\\displaystyle
={\phi_1(x)\over \phi_2(x)}\left({\phi_2(x)\over \phi_1(x)}\right)'\cdot o\left({\phi_2(x)\over \phi_1(x)}\right)=o\left(\left({\phi_2(x)\over \phi_1(x)}\right)'\right),\  x\to x_0.\end{cases}\leqno(9.31)$$

Relation (4.29)$_2$ follows at once from (9.30), (9.31).\ (II)\ The above calculations are valid until (9.30), and in  (9.31) ``$o$'' is replaced by ``$O$'' so that we can only infer (9.19). Now, under the additional conditions, we put
$$
h(x):={f(x)\over\phi_1(x)}-a_1-a_2\phi_(x),\ \ where\ \ \phi(x):=\phi_2(x)/\phi_1(x),\leqno(9.32)$$
so that (9.1) is equivalently written as
$$
h(x)=o\big(\phi(x)\big),\  x\to x_0,\leqno(9.33)$$
and (9.30) yields
$$
h'(x)\equiv \left({f(x)\over\phi_1(x)}\right)'-a_2\phi'(x)=O\big(\phi'(x)\big),\  x\to x_0.\leqno(9.34)$$

To prove $h'(x)=o\big(\phi'(x)\big)$ we try to use some interpolatory-type result on formal differentiation by evaluating $h''$ by a suitable device. Replacing $\phi'=W\cdot (\phi_1)^{-2}$ in factorization (5.10) we get
$$
L[f]={W\over\phi_1}\left[{1\over\phi'}\left({f\over\phi_1}\right)'\right]'=
{W\over\phi_1}\left[\left({1\over\phi'}\right)'\left({f\over\phi_1}\right)'+
{1\over\phi'}\left({f\over\phi_1}\right)''\right],\leqno(9.35)$$
whence
$$\begin{cases}\displaystyle
\left({f\over\phi_1}\right)''=\phi'\left[{\phi_1\over W}L[f]-\left({1\over\phi'}\right)'\left({f\over\phi_1}\right)'\right]=
{\phi_1\phi'\over W}L[f]+{\phi''\over\phi'}\left({f\over\phi_1}\right)' =\\ \\\displaystyle
\overset{\ by\ (9.17)\ and\ (9.19)\ }{=}O\left({\phi' W\over \phi_1\phi_2}\right)+O(\phi'')=O\left({(\phi')^2\over\phi}\right)
+O(\phi'')\overset{\ by\ (9.23)\ }{=}O(\phi''),
\end{cases}\leqno(9.36)$$
and
$$h''\equiv\left({f\over\phi_1}\right)''-a_2\phi''=O(\phi'').\leqno(9.37)$$

Hence $h$ satisfies ``$h=o(\phi),\ h''=O(\phi'')$''. The assumptions on $\phi$ make applicable a result by Boas [15; th. 1B, p. 638] so inferring
$$
h'(x)=o\big(|\phi(x)\phi''(x)|^{1/2}\big)\overset{\ by\ (9.23)\ }{=}o\big(\phi'(x)\big)
,\  x\to x_0.\leqno(9.38)$$
\hfill $\Box$

The original proof  given by Boas is for a function $h$ of class $C^2$ but it can be easily adapted to the case $h\in AC^1$.

As far as (9.21) is concerned we point out that the meaning of  a condition like 
$$
x\phi'(x)\asymp\phi(x),\ x\to+\infty,\  x\to x_0,\leqno(9.39)$$
is properly understood in the context of regular variation. Referring to the monograph by Bingham, Goldie and Teugels [14; th. 2.2.6, p. 74] the positive absolutely-continuous functions $\phi$ satisfying (9.39) are a subclass of the so-called ``{\it extended regularly-varying functions at }$+\infty$'' in the sense of Karamata. A still narrower class is obviously that of functions $\phi$ such that 
$$\phi'(x)/\phi(x)=\alpha x^{-1}+o(x^{-1}),\ x\to +\infty,\ for\ some\ \alpha\in{\mathbb R},\leqno(9.40)$$
 which may be called ``{\it regularly-varying functions at $+\infty$ in a strong sense with index $\alpha$}''.

The analogous class at a point $x_0^-$, where $x_0\in{\mathbb R}$, is defined requiring that the associated functions $\overline{\phi}(x):=
\phi\big((x_0-x)^{-1}\big)$ satisfies (9.40). For such classes of functions a nice result states the equivalence between the pair (4.29)$_{1,2}$ and the pair (9.1)-(9.2). We need an intermediary result.

\vspace{10pt}
{\bf Proposition 9.2.}\  {\it Let hypotheses $(4.13)$ hold true; under the additional assumptions:
$$
\phi_2'(x)=o\big(\phi_1'(x)\big),\ x\to x_0,\leqno(9.41)$$
$$
\phi_1'(x)/\phi_1(x)=O\big(\phi_2'(x)/\phi_2(x)\big),\ x\to x_0,\leqno(9.42)$$
 the pair $(4.29)_{1,2}$ implies the pair} (9.1)-(9.2). 

\vspace{5pt}
{\it Proof.}\ From (4.29)$_2$ we get:
$$
f'\phi_1=f\phi_1'+a_2[\phi_2'\phi_1-\phi_2\phi_1']+o(\phi_2'\phi_1-\phi_2\phi_1')
\overset{\ (4.29)_1\ }{=}
a_1\phi_1\phi_1'+a_2\phi_2'\phi_1+o(\phi_2'\phi_1)+o(\phi_2\phi_1'),$$
from whence
$$
f'=a_1\phi_1'+a_2\phi_2'+o(\phi_2')+o(\phi_2\phi_1'/\phi_1)\overset{\ (9.42)\ }{=}
a_1\phi_1'+a_2\phi_2'+o(\phi_2').\leqno(9.43)$$
\hfill $\Box$

\vspace{10pt}
{\bf Theorem 9.3.}\ {\it Let hypotheses $(4.13)$ and $(9.41)$ hold true together with the following additional assumption:
$$\begin{cases}
\phi_1,\phi_2\ regularly-varying\ at\ x_0^- in\ a\ strong\ sense\\ with\ respective\ indexes\ \alpha_1,\alpha_2;\ \alpha_1>\alpha_2,\ \alpha_2\ne0.
\end{cases}\leqno(9.44)$$

Then  the pair $(4.29)_{1,2}$ is equivalent to the pair} (9.1)-(9.2). 

\vspace{5pt}
{\it Proof.}\ Case: $x_0=+\infty$.\ From (9.40) applied to $\phi_1,\phi_2$  we easily infer the following relations:
$$
\big(\phi_2(x)/\phi_1(x)\big)'\sim (\alpha_2-\alpha_1)x^{-1}\phi_2(x)/\phi_1(x),\ x\to+\infty;\leqno(9.45)$$
$$
\phi_1'(x)/\phi_1(x)\sim  (\alpha_1/\alpha_2)\phi_2'(x)/\phi_2(x),\ x\to+\infty,\ \ if\ \ 
 \alpha_1,\alpha_2\ne0;\leqno(9.46)$$
$$
\phi_1'(x)/\phi_1(x)=o\big(\phi_2'(x)/\phi_2(x)\big),\ x\to+\infty,\ \ if\ \  \alpha_1=0,\alpha_2\ne0.\leqno(9.47)$$

Now the inference ``$(4.29)_{1,2}\implies$ (9.1)-(9.2)'' follows from Proposition 9.2 as
condition  (9.42) is implied by (9.46)-(9.47). Viceversa from (9.1)-(9.2) we get:
$$
f'\phi_1-f\phi_1'=a_2[\phi_2'\phi_1-\phi_2\phi_1']+o(\phi_1\phi_2')+o(\phi_1'\phi_2),
\leqno(9.48)$$
from whence
$$
(f/\phi_1)'=a_2(\phi_2/\phi_1)'+o(\phi_2'/\phi_1)+o\big(\phi_1'\phi_2/(\phi_1)^2
\big).\leqno(9.49)$$

In the case $\alpha_1,\alpha_2\ne0$ we get:
$$
\phi_2'/\phi_1\overset{\ (9.46)\ }{\sim}\displaystyle {\alpha_2\over\alpha_1}\phi_1'\phi_2/(\phi_1)^2\overset{\ (9.44)\ }{\sim}
\alpha_2x^{-1}\phi_2/\phi_1\overset{\ (9.45)\ }{\sim}{\alpha_2\over\alpha_2-\alpha_1}
(\phi_2/\phi_1)',\leqno(9.50)$$
and (4.29)$_2$ follows from (9.49) and (9.450). In the case $\alpha_1=0$ we have:
$$\begin{cases}
\phi_1'\phi_2/(\phi_1)^2\overset{\ (9.47)\ }{=}o(\phi_2'/\phi_1),\\ \\
\phi_2'/\phi_1\overset{\ (9.44)\ }{\sim} \alpha_2 x^{-1}\phi_2/\phi_1
\overset{\ (9.45)\ }{\sim} \displaystyle{\alpha_2\over\alpha_2-\alpha_1}
(\phi_2/\phi_1)',\end{cases}\leqno(9.51)$$
and (4.29)$_2$ follows from (9.49) and (9.51). The case $x_0\in{\mathbb R}$ reduces to the case $x_0=+\infty$ by the mentioned change of variable: asymptotic relations as $x\to x_0^-$ change into equivalent asymptotic relations as $x\to +\infty$, at least in our present situation involving only first-order derivatives. \hfill $\Box$

\vspace{5pt} {\it Remark.}\ The additional conditions in the last two propositions are merely sufficient for the respective theses.
In the elementary case of the scale ``$x\gg1,\ x\to+\infty$'', it happens that the pair $(4.29)_{1,2}$ implies the pair (9.1)-(9.2), and even the stronger relation 
$f'=a_1+o(x^{-1})$, though (9.42) does not hold.

\par\vspace{15pt}
\centerline{\bf References}
\begin{description}
    \item[{[1]}]
    {\small W. A. C}{\scriptsize OPPEL,} {\small \emph{Disconjugacy.} Lecture Notes in Mathematics, vol. 220. 
    Springer-Verlag, Berlin, 1971.}
    \item[{[2]}]
    {\small J. D}{\scriptsize IEUDONNE$^{^{\acute{}}}$},
    {\small \emph{Calcul infinit\'esimal.} Hermann, Paris, 1968.}
    \item[{[3]}]
    {\small A. G}{\scriptsize RANATA,} {\small Canonical factorizations of disconjugate
    differential operators, \emph{SIAM J. Math. Anal}., \textbf{11} (1980),
    160-172.}
    \item[{[4]}]
    {\small A. G}{\scriptsize RANATA,} {\small Canonical factorizations of disconjugate
    differential operators -Part II, \emph{SIAM J. Math. Anal}, \textbf{19} (1988),
    1162-1173.}
    \item[{[5]}]
    {\small A. G}{\scriptsize RANATA,} {\small Polynomial asymptotic expansions in the real domain:
    the geometric, the factorizational, and the stabilization approaches,
    \emph{Analysis Mathematica}, \textbf{33} (2007), 161-198.}
    \item[{[6]}]
    {\small A. G}{\scriptsize RANATA,} {\small The problem of differentiating an asymptotic expansion in real powers. Part I: Unsatisfactory or partial     results by classical approaches,
    \emph{Analysis Mathematica}, \textbf{36} (2010), 85-112.}
        \item[{[7]}]
    {\small A. G}{\scriptsize RANATA,} {\small The problem of differentiating an asymptotic expansion in real powers. Part II: 
    factorizational theory,
    \emph{Analysis Mathematica}, \textbf{36} (2010), 173-218.}
    \item[{[8]}]
    {\small O. H}{\scriptsize AUPT,} {\small \"Uber  Asymptoten ebener Kurven, \emph{J. Reine Angew. Math.},           
    \textbf{152} (1922), 6-10; ibidem p.239.}
    \item[{[9]}]
    {\small E. H}{\scriptsize EWITT} {\small and K. S}{\scriptsize TROMBERG,}
    {\small \emph{Real and abstract Analysis}.
    Springer-Verlag, Berlin, Heidelberg, New York, 1969.}
    \item[{[10]}]
    {\small S. K}{\scriptsize ARLIN} {\small and W. S}{\scriptsize TUDDEN,}
    {\small \emph{Tchebycheff systems: with
    applications in analysis and statistics}. Interscience, New
    York, 1966.}
    \item[{[11]}]
    {\small A. M. O}{\scriptsize STROWSKI,}
    {\small Note on the Bernoulli-L'Hospital rule, \emph{Amer. Math. Monthly}, \textbf{83} (1976), 239-242.}
    \item[{[12]}]
     {\small A. M. O}{\scriptsize STROWSKI,}
    {\small On Cauchy-Frullani integrals, \emph{Comm. Math. Helv.}, \textbf{51} (1976), 57-91.}
    \item[{[13]}]
    {\small W. F. T}{\scriptsize RENCH.,} {\small Canonical forms and principal systems for
    general disconjugate equations, \emph{Trans. Amer. Math. Soc}., \textbf{189}
    (1974), 139-327.}
\end{description}
 
 {\bf Additional references for \S9}
\begin{description}
\item[{[14]}]
{\small N. H. B}{\scriptsize INGHAM,}\ {\small C. M. G}{\scriptsize OLDIE,} and 
{\small J. L. T}{\scriptsize EUGELS,} \emph{Regular variation}. {\small Cambridge University Press, Cambridge, 1987.}
\item[{[15]}]
{\small R. P. B}{\scriptsize OAS Jr.,} {\small Asymptotic relations for derivatives,
\emph{Duke Math. J.,} {\bf 3} (1937), 637-646.}
\item[{[16]}]
{\small D. G. S}{\scriptsize AARI,} {\small An elementary Tauberian theorem for absolutely continuous functions and for series, SIAM \emph{J. Math. Anal.,} {\bf (5)} (1974), 649-662.}
\end{description}

\newpage
{\bf Corrections of misprints in the published version of the present e-paper:}
\\ \\
  {\small A. G}{\scriptsize RANATA,} {\small Analytic theory of finite asymptotic expansions in the real domain. Part I: two-term expansions of differentiable functions, \emph{Analysis Mathematica}, \textbf{37}(2011), 245-287.}
  DOI: 10.1007/s10476-011-0402-7.

  \vspace{10pt}  
On p. 251 in the unnumbered formula $f_i^{\ast}(t)=q_i(t)\cdot L[f(t)]$ correct the left-hand side into  $(f_i^{\ast})'(t)$.

On p. 255 in formula (4.4) the quantity $(f(t)/\phi_1(t_0))'$ must be read as $(f(t)/\phi_1(t))'$.

On p. 259 notice that relation (4.27) in Theorem 4.5 is a different formulation of relation (4.29)$_2$.

On p. 261, in the first line locution \textit{"t limit position"} must be simply read  \textit{"limit position"}.

On p. 263, hypothesis (ii) in Lemma 5.1 must be read "$f\in AC^1(I)$" as stated at the outset of \S5.

On p. 266, in each of the two formulas (5.27)-(5.28) there is a redundand 'absolute value' between the functions inside the integral.

On p. 284, inside the second integral in formula (8.26) the number "$-1$" must be changed into "$+1$", hence the correct version of formula (8.26) is:
$$
\int\limits_{\to 0}^1 t^{1-\alpha_1}\cdot L_{\alpha_1,\alpha_2} [f(t)] dt\equiv \int\limits_{\to 0}^1 [t^{\alpha_2-\alpha_1+1}(t^{-\alpha_2} f(t))']' dt\;\, converges.$$

On p. 286: in reference [7] the page numbers of the paper are missing, namely "173-218".

\end{document}